\documentclass[10pt, letterpaper]{amsart}

\usepackage{amsmath,amssymb,amsfonts,amsthm,mathrsfs, leftidx}
\usepackage{enumerate}
\usepackage[all]{xy}
\usepackage{mathpazo,color}
\usepackage[usenames,dvipsnames]{xcolor}
\usepackage[paper = letterpaper, left = 3cm, right = 3cm, headsep = 6mm, footskip = 10mm,
top = 3cm, bottom = 3cm, footnotesep=5mm, headheight = 2cm]{geometry}
\usepackage{fancyhdr}

\usepackage{comment}

\usepackage[
            breaklinks=true,
            pdftitle={},
            pdfauthor={}]{hyperref}
\usepackage{color}
\definecolor{tocolor}{rgb}{.1,.1,.5}
\definecolor{urlcolor}{rgb}{.2,.2,.6}
\definecolor{linkcolor}{rgb}{.1,.4,.6}
\definecolor{citecolor}{rgb}{.6,.3,.1}
\hypersetup{colorlinks=true, urlcolor=urlcolor, linkcolor=linkcolor, citecolor=citecolor}

\newcommand{\Z}{\mathbb{Z}}

\newcommand{\N}{\mathbb{N}}
\newcommand{\C}{\mathbb{C}}
\newcommand{\R}{\mathbb{R}}
\renewcommand{\O}{\mathcal{O}}
\newcommand{\K}{\mathcal{K}}

\newcommand{\MM}{\mathcal{M}}
\newcommand{\A}{\mathcal{A}}
\renewcommand{\O}{\mathcal{O}}
\renewcommand{\P}{\mathbb{P}}
\newcommand{\B}{\textnormal{B}}

\newcommand{\ii}{\mathbf{i}}
\newcommand{\kk}{\mathbf{k}}
\newcommand{\delbar}{\bar{\partial}}
\newcommand{\dbar}{\bar{\partial}}
\newcommand{\dR}{\textnormal{dR}}

\renewcommand{\a}{\mathfrak{a}}
\newcommand{\g}{\mathfrak{g}}

\newcommand{\p}{\mathfrak{p}}

\DeclareMathOperator{\Spec}{Spec}

\DeclareMathOperator{\End}{End}
\DeclareMathOperator{\Ext}{Ext}

\DeclareMathOperator{\coker}{coker}

\DeclareMathOperator{\Aut}{Aut}
\DeclareMathOperator{\holIso}{Iso_{\mathcal{O}}}
\DeclareMathOperator{\Pic}{Pic}

\DeclareMathOperator{\Sym}{Sym}

\newcommand{\Ad}{\textnormal{Ad}}

\newcommand{\arsim}{\xrightarrow{\sim}}
\newcommand{\codim}{\textnormal{codim}}
\newcommand{\commsq}[8]{\xymatrix{ #1  \ar[r]^{#5} \ar[d]_{#6} & #2 \ar[d]^{#7} \\ #3 \ar[r]_{#8} & #4 }}

\newcommand{\mat}[4]{\left[ \begin{array}{cc} #1 & #2 \\ #3 & #4 \end{array} \right]}

\newcommand{\acts}{\mbox{\,\raisebox{0.26ex}{\tiny{$\bullet$}}\,}}

\theoremstyle{plain}
\newtheorem{thm}[equation]{Theorem}
\newtheorem{lem}[equation]{Lemma}

\newtheorem{question}[equation]{Question}
\newtheorem{prop}[equation]{Proposition}
\newtheorem{cor}[equation]{Corollary}

\theoremstyle{definition}
\newtheorem{defn}[equation]{Definition}
\newtheorem{ex}[equation]{Example}

\theoremstyle{remark}
\newtheorem{rmk}[equation]{Remark}

\numberwithin{equation}{section}

\begin{document}

\title{Canonical complex extensions of K\"ahler manifolds}

\author{Daniel Greb}
\address{Essener Seminar f\"ur Algebraische Geometrie und Arithmetik, Fakult\"at f\"ur Mathematik, Universit\"at Duisburg--Essen, 45117 Essen, Germany}
\email{daniel.greb@uni-due.de, michael.wong@uni-due.de}
\author{Michael Lennox Wong}

\date{\today}

\begin{abstract} 
Given a complex manifold $X$, any K\"ahler class defines an affine bundle over $X$, and any K\"ahler form in the given class defines a totally real embedding of $X$ into this affine bundle. We formulate conditions under which the affine bundles arising this way are Stein and relate this question to other natural positivity conditions on the tangent bundle of $X$. For compact K\"ahler manifolds of non-negative holomorphic bisectional curvature, we establish a close relation of this construction to adapted complex structures in the sense of Lempert--Sz\H{o}ke and to the existence question for good complexifications in the sense of Totaro. Moreover, we study projective manifolds for which the induced affine bundle is not just Stein but affine and prove that these must have big tangent bundle. In the course of our investigation, we also obtain a simpler proof of a result of Yang on manifolds having non-negative holomorphic bisectional curvature and big tangent bundle.
\end{abstract}

\maketitle

\section*{\textbf{Introduction}}

Given a K\"ahler class $[\omega] \in H^{1,1}\bigl(X \bigr) \cong H^1 \bigl( X, \, \Omega^1_X\bigr) \cong \mathrm{Ext}^1\bigl(\mathcal{O}_X,\, \Omega_X^1 \bigr)$ on a complex manifold $X$, the last identification allows one to construct a natural deformation of the complex structure on the cotangent bundle of $X$, which is an affine bundle modelled on $\Omega_X^1$, see for example \cite[\S 2]{Donaldson2002}. In this paper, we study the complex geometry of the total space of this affine bundle $Z_X = Z_{[\omega]}$ over $X$, which we call a "canonical complex extension of $X$", adopting the terminology of \cite{S} and referring the reader to the discussion in Section~\ref{s:cxfn} for why we have chosen the term "extension" over the perhaps more natural "complexification".

It is elementary to see that every form representing $[\omega]$ defines a differentiable section of $Z_{[\omega]}$. We show in Theorem~\ref{thm:totallyreal} below that a section obtained in this way defines a totally real embedding of $X$ into $Z_{[\omega]}$ if the chosen form is K\"ahler. As every K\"ahler class has a real analytic K\"ahler representative, a folklore result which we prove in detail in Appendix~\ref{app:anrep} of this paper, such a section hence provides a natural complex extension of the real analytic manifold underlying $X$. In particular, by classical results in complex analysis due to Grauert, the image of such a section has a neighbourhood basis consisting of Stein domains, and it is hence a natural question to ask under which conditions is $Z_{[\omega]}$ itself Stein. 

Looking at what happens in the case of compact Riemann surfaces is perhaps useful to form an initial set of expectations.  For $X = \P^1$, $Z_X$ is an affine variety (in fact, isomorphic to the surface in $\C^3$ defined by $x^2 + y^2 + z^2 = 1$, hence a complexification of $S^2 \cong \P^1$), and when $X$ is an elliptic curve, $Z_X$ is Serre's example of an algebraic variety which is not affine, but whose underlying complex structure is Stein, see Remark \ref{r:notalgbutStein}.  Finally, when the genus of $X$ is greater than or equal to two, $Z_X$ admits no non-constant holomorphic functions, and hence is never Stein, see Example \ref{ex:RS}.  On the other hand, $Z_X$ (for $X$ of arbitrary dimension) always possesses the property of Stein manifolds that there are no compact complex submanifolds of positive dimension, see Proposition \ref{p:0dim}.

\subsection*{Steinness for manifolds with nonnegative holomorphic bisectional curvature}
Based on the foregoing ---admittedly crude---set of data, it is plausible to suggest that some positivity conditions on the tangent bundle of $X$ may be sufficient for $Z_X$ to be a Stein manifold.  After studying the examples of complex tori and flag manifolds in Sections~\ref{s:Euctori}, we prove that for K\"ahler manifolds $(X, \omega)$ admitting a finite \'etale covering fibering over a compact complex torus with fibre a flag manifold, the affine bundle $Z_{[\omega]}$ is indeed Stein, see Theorem~\ref{t:bisholStein}.  As a consequence of Mok's uniformisation theorem \cite[Main Theorem, p.179]{Mok1988}, this class of manifolds includes all compact K\"ahler manifolds of non-negative holomorphic bisectional curvature, see Proposition \ref{p:bisfinitecover}.

\subsection*{Relation to adapted complex structures}
There is another canonical way of associating to a Riemannian metric $g$ on a given real analytic manifold $M$ a complex structure on a neighbourhood of the zero section of the tangent bundle $TM$, the so-called \emph{adapted complex structure} of Lempert--Sz\H{o}ke~\cite{LempertSzoke1991} and Guillemin--Stenzel \cite{GuilleminStenzel1991}. In fact, one of the natural questions in this theory is to determine conditions under which the adapted complex structure is everywhere defined, so that it induces the structure of a Stein manifold on the entirety of the tangent bundle.  When this happens, the Riemannian manifold $(M, g)$ is said to have an \emph{entire Grauert tube}, and one knows that it has non-negative sectional curvature \cite[Theorem 2.4]{LempertSzoke1991}.  On the other hand, it is not true that non-negative sectional curvature of a particular metric is sufficient for the complex structure arising from that metric to be defined on the whole of $TM$.\footnote{It is believed, however, that if there exists a metric of non-negative sectional curvature, then there is some, possibly different, metric, necessarily also of non-negative curvature, for which the Grauert tube is entire; this appears to be a subtle and difficult question, related to the suggestion of Totaro mentioned below.}  Of course, in the case that a Riemannian metric of non-negative sectional curvature is also K\"ahler, then its holomorphic bisectional curvature is also non-negative.  In this K\"ahler situation, we will see that one often already has an entire Grauert tube, and even if this is not guaranteed, there is always another metric in the same K\"ahler class for which this holds; furthermore, this adapted complex structure is then biholomorphic to the canonical complex extension $Z_{[\omega]}$.  See Theorem~\ref{thm:same_as_adapted} for a precise statement.

\subsection*{Good complexifications}
In fact, a conjecture of Burns \cite[\S5]{Burns1982} predicts that for a Riemannian manifold $(M, g)$ with entire Grauert tube, the tangent bundle can be made into an affine algebraic variety in a natural way. Already for tori this is a subtle question, as for example our computations of the complex structure obtained by deformation on tangent bundles of compact complex tori in Section~\ref{s:Euctori} below show that in this context Serre's classical examples of non-affine algebraic varieties with Stein analytification do show up. This algebraisation problem is closely related to the existence of a \emph{good complexification} of the underlying real manifold $M$, i.e., the question of whether there exists a smooth affine algebraic variety $Z$ over $\mathbb{R}$ such that $M$ is diffeomorphic to $Z(\mathbb{R})$ and such that $Z(\mathbb{R}) \hookrightarrow Z(\mathbb{C})$ is a homotopy equivalence.

Totaro suggests that manifolds admitting a Riemannian metric of non-negative sectional curvature should be precisely those that admit such a good complexification, see \cite[Introduction]{Totaro2003} for a detailed discussion.  Using some of the techniques introduced in \cite{Totaro2003} as well as the rough structure theory provided by Proposition \ref{p:bisfinitecover}, we show that in the K\"ahler case non-negative curvature (hence, non-negative holomorphic bisectional curvature) is indeed sufficient for the existence of good complexifications, see Theorem~\ref{thm:good_complexification}. 

\subsection*{Necessary criteria}
Having discussed sufficient criteria for Steinness of the affine bundle $Z_{[\omega]}$, it is of course an interesting question to characterise those compact K\"ahler manifolds fulfilling this positivity condition.  In this direction, note that compact K\"ahler manifolds with non-negative holomorphic bisectional curvature constitute an important class of compact K\"ahler manifolds with nef tangent bundle, cf.\ \cite{DPS}. As the tautological line bundle on the projectivisation of the tangent bundle appears as the normal bundle of a divisor compactifying $Z_{[\omega]}$ to a compact K\"ahler manifold, see Lemma~\ref{l:extlem}, it is natural to ask whether compact K\"ahler manifolds for which $Z_{[\omega]}$ is Stein do in fact have nef tangent bundle. Since proving this seems to be difficult, as a first step in this direction, we investigate the algebraic situation and show in Section~\ref{subsect:big_tangent} that for projective manifolds having a class ${[\omega]}$ such that $Z_{[\omega]}$ is affine, the tangent bundle is necessarily big. As a side product of our considerations we obtain essentially elementary proofs of a result of Hsiao \cite[Corollary 1.3]{Hsiao2015}, asserting that flag varieties have big tangent bundles, and one of Yang \cite[Theorem 4.5]{Yang2017}, stating that a compact K\"ahler manifold of non-negative holomorphic bisectional curvature with a big tangent bundle is a product of irreducible Hermitian symmetric spaces of compact type and of projective spaces with metrics of non-negative holomorphic bisectional curvature.

\subsection*{Future directions}
Manifolds with nef tangent bundle necessarily have nef anti-canonical bundle; it therefore follows from \cite[Cor.~1.4]{Cao} and \cite[Main Theorem]{DPS} that the universal cover of such a manifold is a direct product of $\mathbb{C}^k$ with a Fano manifold with nef tangent bundle.  It was conjectured by Campana--Peternell in \cite{CP} that such Fano manifolds are in fact rational homogeneous projective manifolds, i.e., flag varieties. Comparing with the proof of Steinness of the canonical extension in the case of manifolds with non-negative bisectional curvature, one is led to the following natural question. 
\begin{question}
 If $X$ is a compact K\"ahler manifold with nef tangent bundle, is there a K\"ahler metric $\omega$ on $X$ such that $Z_{[\omega]}$ is Stein? 
\end{question}

While the discussion above shows that the universal covers of manifolds with non-negative bisectional curvature and of those with nef tangent bundle modulo the Campana--Peternell conjecture are the same as complex manifolds, we remark that the statement of Proposition~\ref{p:bisfinitecover} does not generalise to the more general setup of manifolds with nef tangent bundle.\footnote{Indeed, the projective bundle $X = \mathbb{P} W$ associated with the non-split extension $0 \to \mathcal{O}_E \to W \to \mathcal{O}_E \to 0$ on an elliptic curve $E$, on the one hand, has nef tangent bundle by \cite[Thm.~3.1]{CP}, while on the other hand by \cite[Rem.~1.7]{CDP}, it does not admit a K\"ahler metric such that the universal cover splits isometrically and biholomorphically as $\mathbb{P}^1 \times \mathbb{C}$, where both factors are endowed with their standard metric.} Concerning the converse problem of characterising manifolds admitting a metric with Stein canonical extension, one might ask

\begin{question}
 Let $X$ be a compact K\"ahler manifolds having a K\"ahler metric $\omega$ such that $Z_{[\omega]}$ is Stein. Does this imply that $\mathcal{O}_{\mathbb{P}\Omega_X^1} (1)$ is pseudoeffective / modified nef in the sense of \cite[Def.2.2]{Boucksom} / nef ?
\end{question}

Finally, while we present no results in this direction, we will end our introduction by mentioning that our considerations in this paper should be related to the hyperk\"ahler geometry of cotangent bundles of compact K\"ahler manifolds, as developed in \cite{Feix, Kaledin2001}.  One of the features of this theory resembles that of the adapted complex structures, in that the hyperk\"ahler metric does not always exist on the whole cotangent bundle, but only in some neighbourhood of the zero section.  Furthermore, the situation for compact Riemann surfaces mirrors that described above:  in genera $0$ and $1$, the hyperk\"ahler metric exists everywhere, while for genus greater than or equal to two, it is known that it does not extend to the whole cotangent bundle \cite[Theorem B(ii)]{Feix}; in fact, in the latter case, the domain of definition of the metric has more recently been determined \cite[\S6.1]{Hitchin2015}.  It is likely that the assumption of non-negative holomorphic bisectional curvature is again sufficient to guarantee that the hyperk\"ahler metric exists on the whole cotangent bundle and is complete, but we have not investigated this.

\subsection*{Acknowledgements}
The authors would like to thank Bo Berndtsson, Andreas H\"oring, and Stefan Nemirovski for answering questions via eMail, Indranil Biswas for pointing out Example \ref{ex:RS}, Tim Kirschner for in-depth discussions regarding the material covered in Appendix~\ref{app:anrep}, L\'aszl\'o Lempert for bringing our attention to \cite{Semmes} and for sharing his related work with us, and Robert Sz\H{o}ke for insightful explanations regarding adapted complex structures and Burns' conjecture, as well as for pointing out an initial misinterpretation on our part of Mok's theorem. Both authors were supported by the Collaborative Research Center SFB/TR 45 ``Periods, moduli spaces and arithmetic of algebraic varieties''(Project M08-10) of the Deutsche Forschungsgemeinschaft (DFG).

\section{\textbf{Preliminaries}}

If not mentioned otherwise, manifolds are assumed to be connected and differentiable manifolds are assumed to be $C^\infty$.

\subsection{Complexification} \label{s:cxfn}

By Whitney's Theorem, every $C^\infty$-manifold admits a real-analytic structure. The complexification of such a real analytic manifold seems to have first appeared in \cite{WhitneyBruhat}.  The arguments in that paper are somewhat sparse and more details are given in \cite[\S1.4]{S}. While these early papers have a local approach to the notion of complexification, we will use the following definition, cf.~\cite{Kulkarni1978}.

\begin{defn} \label{d:cxfn}
 A \textbf{complexification} of a real analytic manifold $M$ is a complex manifold $Z$ together with an anti-holomorphic involution $\tau : Z \to Z$ such that $M$ is real analytically isomorphic to the fixed-point set $Z^\tau$ of $\tau$.  The complexification $Z$ is called \textbf{minimal} if the inclusion $M \hookrightarrow Z$ is a homotopy equivalence.
\end{defn}
\begin{rmk}
 Two complexifications of the same real analytic manifold $M$ have biholomorphic germs around $M$.
\end{rmk}

\begin{defn} 
Let $M$ be a differentiable manifold and $Z$ a complex manifold.  We say that an embedding $\imath : M \hookrightarrow Z$ is \textbf{totally real} if for all $x \in M$, $d\imath_x(T_x M) \cap J_{Z, i(x)} \cdot d\imath_x(T_x M) = 0$, where $J_Z : TZ \to TZ$ is the complex structure on the tangent bundle. In this case, we call $Z$ a \textbf{(complex) extension} of $M$.
\end{defn}
Fixed-point sets of anti-holomorphic involutions are totally real, i.e., a complexification is a complex extension. Conversely, Grauert's solution of the Levi problem \cite{Grauert1958} and further improvements by Nirenberg--Wells \cite[\S 4]{NirenbergWells} imply the following result. 

\begin{thm}[Real analytic totally real embeddings and complexifications] \label{thm:GrauertNirenberg}
Every real analytic totally real embedding $M\hookrightarrow Z$ admits a neighbourhood basis of Stein open subsets $M \subset \Omega_\varepsilon \subset Z$ such that each $\Omega_\varepsilon$ is a minimal complexification of $M$. 
\end{thm} 
In particular, every complexification contains a minimal complexification that is a Stein manifold. As smooth affine varieties over the complex numbers provide examples of Stein manifolds, this leads to the following notion, which was introduced in \cite{Totaro2003}, and which we discuss in our context in Section \ref{s:goodcxfn}. 

\begin{defn}
 A \textbf{good complexification} of a real analytic manifold $M$ is a smooth affine algebraic variety $U$ over $\R$ such that $M$ is diffeomorphic to $U(\R)$ and such that the inclusion $U(\R) \to U(\C)$ is a homotopy equivalence. 
\end{defn}

We will later use the following uniqueness property of complex extensions:

\begin{lem} \label{l:uniquelift}
Let $\psi : M \to M$ be a diffeomorphism of the differentiable manifold $M$, and let $\imath : M \hookrightarrow Z$ be a complex extension of $M$. Then, there is at most one holomorphic automorphism $\Psi : Z \to Z$ making $\imath$ equivariant with respect to $\psi$ and $\Psi$.   In particular, if a finite group $\Gamma$ acts by diffeomorphisms on $M$, then there is at most one holomorphic action of $\Gamma$ on $Z$ such that $\imath$ is $\Gamma$-equivariant.
\end{lem}

\begin{proof}
It suffices to show that a holomorphic automorphism of $Z$ that fixes $i(M)$ pointwise is the identity. As $i(M)$ is totally real inside $Z$, the Cauchy--Riemann equations imply that such an automorphism is the identity on a neighbourhood of $i(M)$ in $Z$. We then conclude using the identity principle for holomorphic maps between (connected) complex manifolds. 
\end{proof}

\subsection{Analytic representatives for K\"ahler forms}

While many papers on complexifications in the presence of a K\"ahler form assume the form to be real analytic, see for example \cite{Feix}, the following result does not seem to be too well-known.  As only sketches of proofs of Proposition~\ref{t:anrep} can be found in the literature, we give a detailed argument in Appendix \ref{app:anrep}.

\begin{prop} \label{t:anrep}
Let $X$ be a K\"ahler manifold with K\"ahler form $\omega_0$.  Then there exists $\psi \in C^\infty(X)$ such that the form $\omega := \omega_0 + i \partial \dbar \psi$ is positive and real analytic.
\end{prop}

\begin{rmk}
Observe that $\omega = \omega_0 - \dbar( i \partial \psi)$, so $\omega$, $\omega_0$ represent the same Dolbeault cohomology class in $H^{1,1}_{\bar\partial}(X)$.  Also, since $i \partial \dbar \psi = d ( \tfrac{i}{2} ( \dbar \psi - \partial \psi))$, they represent the same de Rham cohomology class.  Thus, the statement says that any K\"ahler class (for either of the cohomology theories) has a real analytic representative.
\end{rmk}

\subsection{Affine bundles} \label{s:affbdl}

Let $X$ be a complex manifold and $V$ a holomorphic vector bundle over $X$.  A class $a \in H^1(X, V) = \Ext^1(\O_X, V)$ determines an affine bundle $Z_a \to X$ over $X$ with bundle of translations $V$ as follows.  Associated to $a$ is an extension of vector bundles
\begin{align} \label{e:Vext}
0 \longrightarrow V \longrightarrow W \longrightarrow \O_X \longrightarrow 0.
\end{align}
To recall, $W$ may be realised via Dolbeault representatives as follows.  Let $\alpha \in \mathcal{A}^{0,1}(V)$ denote a Dolbeault representative of $a$.  Then $W$ is the $C^\infty$ vector bundle $V \oplus \O_X$ with the holomorphic structure
\begin{align} \label{e:Dolstr}
\delbar_W := \mat{ \delbar_V }{ \alpha }{}{ \delbar_{\O_X} }.
\end{align}
In terms of \v{C}ech representatives, if $\{ X_i \}$ is a fine enough open cover so that $a$ is represented by $( a_{ij}) \in \prod \Gamma(X_i \cap X_j, V)$, then a section of $W$ over $U$ is given by tuples 
\begin{align} \label{e:WsnCech}
(s_i, f_i) \in \prod_i \Gamma(U \cap X_i, V) \oplus \Gamma(U \cap X_i, \O_X)
\end{align}
satisfying
\begin{align} \label{e:Wsn}
\left[ \begin{array}{c}
s_i \\ f_i \end{array} \right] = \mat{ \mathrm{Id} }{ a_{ij} }{}{ \mathrm{Id} } \left[ \begin{array}{c} s_j \\ f_j \end{array} \right]
\end{align}
on $U \cap X_i \cap X_j$.  Note that if $\alpha$ is the Dolbeault representative taken above, then since $\alpha$ is $\delbar_V$-closed, by the $\delbar$-Poincar\'e lemma, one can find an open cover $\{ X_i \}$ such that there exist $\rho_i \in \A^0(X_i, V)$ with $\delbar_V \rho_i = \alpha|_{X_i}$.  With this, one can take $(a_{ij} = \rho_i - \rho_j)$ for the \v{C}ech representative above.

Let us be explicit about how one goes back and forth between expressions for (not necessarily holomorphic) sections of $W$ in the Dolbeault and in the \v{C}ech realisations.  As mentioned, a section for the \v{C}ech realisation is given by a tuple as in \eqref{e:WsnCech} satisfying \eqref{e:Wsn}, while one of the Dolbeault realisation is a smooth section of the direct sum $V \oplus \O_X$, and this will be holomorphic if and only if it lies in the kernel of \eqref{e:Dolstr}.  Starting with an expression as in \eqref{e:WsnCech}, the corresponding Dolbeault section is given by
\begin{align} \label{e:CechtoDol}
(s_i - \rho_i f_i, f_i)
\end{align}
over $U \cap X_i$; it is easy to check that these expressions give a well-defined section of $V \oplus \O_X$.  Of course, given a Dolbeault section $(s,f)$, the corresponding \v{C}ech section is $(s|_{U \cap X_i} + f|_{U \cap X_i} \rho_i, f|_{U \cap X_i} )$.

If $|W|$ denotes the total space of $W$, then the surjective map in \eqref{e:Vext} yields one of manifolds
\begin{align*}
p : |W| \to |\O_X| = X \times \C \to \C.
\end{align*}
For $\lambda \in \C$, we set $Z_{a, \lambda} := p^{-1}(\lambda)$ and $Z_a := Z_{a, 1}$; of course, we may simply write $Z$ if the class $a$ is understood.  In fact, it is easy to see that $Z_a \cong Z_{a, \lambda}$ for any $\lambda \in \C^\times$.  Also, since the isomorphism class of $W$ only depends on the class of $a$ in $\P H^1(X,V)$, the same is true of $Z_a$.  The vector bundle $V$ acts on $W$ by translations; since $V$ is the kernel of $p$, it is clear that $Z_a$ is invariant (as a submanifold of $|W|$) under this operation, so that we obtain a well-defined action of $V$ on $Z_a$. With this operation, $Z_a$ becomes a $V$-torsor; this is equivalent to saying that $Z_a$ is a locally trivial fibre bundle with fibre $\C^r$ with transition functions in the group of affine transformations of $\C^r$ such that the induced cocycle with values in $GL_r(\C)$ gives back $V$. We say that $Z_a$ is an \textbf{affine bundle (modelled on the vector bundle $V$)}.

\begin{rmk} \label{r:Pconv}
To avoid any possibility of ambiguity, in the following, $\P W$ will refer to the projective bundle of \emph{lines, and not hyperplanes}, in $W$. This is also the convention adopted in \cite[\S 15.C]{DemaillyBook}.
\end{rmk}

The next two lemmata collect a number of well-known properties of the construction described above. We explicitly state and prove them here in order to point out a number of consequences of the choice of convention adopted here, to remind the reader of certain isomorphism between the different realisations, and to be able to easily apply them later in the special case where we consider the extension of the trivial line bundle by the cotangent bundle induced by a K\"ahler class.

\begin{lem}[Sections of $Z_a$] \label{l:affsn}
\begin{enumerate}[(a)]
\item \label{l:cansn} A ($C^\infty$) section of $Z_a$ over $U$ is equivalent to a tuple $(s_i)$ of sections of $V$ over $U \cap X_i$ satisfying the relation
\begin{align*}
s_i = s_j + a_{ij} \quad \text{ on }X_i \cap X_j.
\end{align*}
In particular, taking $s_i = \rho_i$, a choice of local primitives for a Dolbeault representative for $a$ yields a canonical $C^\infty$ section (which will not, in general, be holomorphic).  In the Dolbeault realisation, this is simply the constant section $(0,1)$ of $V \oplus \O_X$.
\item \label{l:sn-triv} A global holomorphic section of $Z_a$ exists if and only if $Z_a \cong V$ if and only if $a \in H^1(X, V)$ is the trivial class.
\item The projection $Z_a \to X$ is a homotopy equivalence.
\end{enumerate}
\end{lem}

\begin{proof}
By definition, a section of $Z_a$ is a section of $W$ whose image in $\O_X$ is $1$.  Thus, the statement comes from taking $f_i = 1$ in \eqref{e:Wsn}.  The statement involving the Dolbeault representation of the canonical section simply comes from \eqref{e:CechtoDol}. Thus statement \eqref{l:cansn} follows.

Regarding (b), we first remark that an affine space is isomorphic to its vector space of translations upon a choice of base point.  Similarly, an affine bundle will be isomorphic to its vector bundle of translations upon choice of a section.  It is clear that if $a$ is trivial, then $V \cong Z_a$.  On the other hand, a global holomorphic section will produce a global primitive for the class of $a$, which will then be the trivial class.

For (c), we noted above that a $C^\infty$ section of $Z_a$ always exists, so $Z_a$ and $V$ are diffeomorphic, but the projection $V \to X$ is a homotopy equivalence, as it is a vector bundle.
\end{proof}

\begin{lem} [Basic properties of the construction] \label{l:extlem}
In the setup adopted in the current subsection, the following holds.
\begin{enumerate}[(a)]
\item \label{l:PWminusPV} The inclusion $V \to W$ in \eqref{e:Vext} induces an inclusion of the bundles $\P V \to \P W$ of projective spaces over $X$.  Then $\P V$ is the vanishing locus of a section $s_0 \in H^0( \P W, \O_{\P W}(1) )$ and  
\begin{align*}
Z_a \cong \P W \setminus \P V.
\end{align*}

\item \label{l:baspropb} Let $f : Y \to X$ be any holomorphic map of complex manifolds and let $V$ be a holomorphic vector bundle over $X$.  Let $a \in H^1(X, V)$ and let $p_a : Z_a \to X$ be the corresponding affine bundle.  The pullback class $f^* a$ lies in $H^1(Y, f^*V)$, and so we also have an affine bundle $p_{f^* a} : Z_{f^* a} \to Y$ and the following is a cartesian diagram:
\begin{align} \label{e:baspropbdiag}
\begin{gathered} \commsq{ Z_{f^* a} }{ Z_a }{ Y }{ X. }{ F }{ p_{f^* a} }{ p_a }{ f } \end{gathered} 
\end{align}

\item \label{l:extlemc} Fix $f : Y \to X$, $V \to X$, and $a \in H^1(X, V)$ as above.  There is a bijection between the set of commutative diagrams
\begin{align} \label{e:extlemcdiag}
\begin{gathered} \xymatrix{ Y \ar[r]^\phi \ar[dr]_f & Z_a \ar[d]^{p_a} \\ & X } \end{gathered}
\end{align}
and the set of sections of $Z_{f^*a} \to Y$ which, by Lemma \ref{l:affsn}\eqref{l:sn-triv} above, is non-empty precisely when $f^* a \in H^1(Y, f^*V)$ is trivial.
\end{enumerate}
\end{lem}

\begin{proof}
To see \eqref{l:PWminusPV}, say in the \v{C}ech realisation, sections of $W$ are locally given by pairs $(s, f)$, with $s$ a section of $V$ and $f$ a locally defined function.  The sections of $V$ correspond to those where $f = 0$.  Therefore, taking $f = 1$ is equivalent to dehomogenising coordinates of $\P W$, and this is precisely the complement of $\P V$.  Recall that under the convention adopted in Remark~\ref{r:Pconv} we have an identification $H^0(\P W, \O_{\P W}(1)) = H^0(X, W^\vee)$; the section $s_0$ is (up to scalars) precisely the section corresponding to the inclusion $\O_X \to W^\vee$ one obtains upon dualising the sequence \eqref{e:Vext}.

For \eqref{l:baspropb}, let $\{ X_i \}$ be an open cover of $X$ for which we can choose a \v{C}ech representative $( a_{ij} ) \in \prod_{i,j} \Gamma(X_i \cap X_j, V)$ for $a$.  Then $Z_a$ is constructed from $\coprod_i V_i$, where $V_i = \big|V|_{X_i}\big|$ is the total space of $V$ over $X_i$, modulo relations built from the $a_{ij}$.

If $Y_i := f^{-1}(X_i)$, then $Z_{f^* a}$ is built similarly from $\coprod W_i$, with $W_i := \big|f^*V|_{Y_i} \big|$ and relations coming from $f^* a_{ij}$.  But virtually by definition, we have $W_i = Y_i \times_{X_i} V_i$. Thus, given a complex manifold $T$ and a commutative diagram
\begin{align*} 
\vcenter{ \commsq{ T }{ Z_a }{ Y }{ X }{ }{ q }{ p_a }{ f } }
\end{align*}
one gets maps $T_i := q^{-1}(Y_i) \to V_i$, and hence to $W_i$, and one will be able to patch these together uniquely into a map $T \to Z_{f^* a}$.

Statement \eqref{l:extlemc} follows directly from the universal property of the cartesian diagram \eqref{e:baspropbdiag}, with a section $s$ of $p_{f^*a}$ corresponding to the map $\phi = F \circ s$ in \eqref{e:extlemcdiag}.
\end{proof}

\begin{rmk}
Lemma \ref{l:extlem}\eqref{l:extlemc} allows us to call $Z_a$ the ``universal space'' on which $a$ trivialises.
\end{rmk}

\subsection{Flag varieties: Automorphisms, isometries, and embeddings} \label{s:fvaut}

By a flag variety, we will mean a smooth complex projective variety of the form $M = G/P$, where $G$ is a connected semisimple complex algebraic group and $P \leq G$ is a parabolic subgroup. 

\subsubsection{Automorphisms} \label{subsubsect:flag_automo}
We wish to remind the reader of the description of the holomorphic/algebraic automorphism group $\Aut M$ of such an $M$; $\Aut^\circ \! M$ will denote the identity component of $\Aut M$.  We follow \cite[\S3.3]{Akhiezer}.

Of course, $G$ itself always acts on $M$ by left multiplication, and as $G$ is connected, we always have a map $G \to \Aut^\circ \! M$.  In the case that $G$ is simple and of adjoint type (that is, with trivial centre), then except for a few special cases \cite[\S3.1. Example 2]{Akhiezer}, this is an isomorphism, so that we may make the identification \cite[\S3.3. Theorem 2]{Akhiezer}
\begin{align} \label{e:Aut0M}
\Aut^\circ \! M = G.
\end{align}
In each of the exceptions, if $G' := \Aut^\circ \! M$, then $G'$ is also a semisimple group, $G \leq G'$ (via the map above) and $G'$ has a parabolic subgroup $P'$ for which $G/P = G'/P'$.  Since we are interested in this latter quotient, by replacing the pair $(G,P)$ with $(G', P')$, we may assume \eqref{e:Aut0M} holds in all cases.

Now, when $G$ is not necessarily simple, one has a decomposition
\begin{align*}
M \cong M_1 \times \ldots \times M_r,
\end{align*}
where each $M_i = G_i/P_i$ is itself a flag manifold, with $G_i$ being a simple factor of $G$ and $P_i = G_i \cap P$, and \cite[\S3.3. Theorem 1]{Akhiezer} says that
\begin{align*}
\Aut^\circ \! M = \Aut^\circ \! M_1 \times \ldots \times \Aut^\circ \! M_r.
\end{align*}
Therefore, by replacing the $G_i$ by the appropriate $G_i'$ as in the preceding paragraph, we may assume that \eqref{e:Aut0M} holds, even if $G$ is not simple.

To describe the full automorphism group, let $T \leq G$ and $B \leq G$ be a maximal torus and a Borel subgroup, respectively, with $T \leq B$; furthermore, we choose them so that $B \leq P$.  Let $\Phi$ denote the root system associated to $(G,T)$; $B$ then corresponds to a subset $\Pi \subseteq \Phi$ of simple roots.  Let $\Aut G$ denote the group of automorphisms of $G$ as an algebraic group and set
\begin{align*}
E := \{ \varphi \in \Aut G \mid \varphi(B) = B, \varphi(T) = T \}.
\end{align*}
Then $E$ is a finite group; since it preserves $T$, it induces an action on the root system $\Phi$; since it furthermore preserves $B$, it permutes the positive roots, and hence is isomorphic to the automorphism group of the Dynkin diagram associated to the \'epinglage of $G$ determined by $T$ and $B$ (i.e., the choice of $\Pi$ above).  Let $E_P \leq E$ be the subgroup of $E$ preserving $P$.  Since $B \leq P$, $P$ corresponds to a subset $R \subseteq \Pi$ of simple roots; viewing $E_P$ as a group acting on the Dynkin diagram, whose nodes are indexed by $\Pi$, it is the subgroup which preserves $R$.

\begin{thm}[Thm.~3 in \S3.3 of \cite{Akhiezer}] \label{t:AutM}
$\Aut M = E_P \ltimes G$.  In particular, $\Aut M$ is an affine algebraic group.
\end{thm}

We will also need the following.

\begin{prop}[Cor.~2 in \S3.3 of \cite{Akhiezer}] \label{p:Hssbas}
Every automorphism of $M$ has at least one fixed point.
\end{prop}

Let $L \leq P$ be the Levi subgroup containing $T$; then, as it is a reductive group itself, its root system $\Phi_L$ with respect to $T$ is a sub-root system of $\Phi$; in fact, it is precisely the sub-root system of $\Phi$ spanned by $R$.  Thus, we obtain the following.

\begin{lem}
If $\varphi \in E_P$, then $\varphi(L) = L$.
\end{lem}

To be explicit, an element of the $G$-factor of $\Aut M = E_P \ltimes G$ acts by left multiplication on $M = G/P$ and an element of $E_P$ acts as an automorphism of $G$ on a representative of a coset in $M$ (of course, it is well-defined as $P$ is preserved by $E_P$, by definition).  The lemma now implies that if we let $Z := G/L$, where $L \leq P$ is the Levi factor described above, then  $E_P \ltimes G$ also acts on $Z$.  Since any algebraic action on $M$ factors through $\Aut M$, we may record the following.

\begin{cor} \label{c:Mactionlift}
Let $H$ be an algebraic group acting algebraically on a flag variety $M = G/P$, where $G = \Aut^\circ \! M$.  Then the $H$-action on $M$ canonically lifts to an algebraic action on $Z$ making the canonical fibration $\pi : Z \to M$ equivariant.
\end{cor}

\subsubsection{Flag varieties in terms of compact groups} \label{s:auts}

A flag manifold may also be realised as a homogeneous space of a compact Lie group.  Let $M = G/P$ be as before (with $G$ semisimple, but not necessarily simple).  Let $U \leq G$ be a maximal compact subgroup (which will also be semisimple); then $G = U^\C$ is the complexification of $U$.  Let $K := U \cap P$.  Then we have identifications
\begin{align*}
M = G/P = U/K.
\end{align*}
Furthermore, we have $K^\C = L$, that is, the complexification of $K$ is the Levi factor of $P$.  With this, the canonical fibration $\pi : Z = G/L \to M$ has a real analytic (even real algebraic), totally real section $s : M \to Z$ given by realising $M = U/K$ and taking
\begin{equation} \label{eq:section_in_GmodL}
 s(uK) := uL.
\end{equation}

Using the description of $\Aut M$ in Theorem \ref{t:AutM} and Lemma \ref{l:uniquelift}, it is straightforward to verify the following.

\begin{lem} \label{l:Meqvtsection}
If an algebraic group $H$ acts algebraically on $M$, then the inclusion $s : M \to Z$ is $H$-equivariant with respect to the lifted action given by Corollary \ref{c:Mactionlift}. Moreover, the latter is uniquely determined by the $H$-action on $M$.
\end{lem}

\subsubsection{Riemannian metrics on flag manifolds}

Given a flag manifold $M = G/P = U/K$, we consider it as a K\"ahler manifold and a fortiori as Riemannian manifold on which $U$ acts by isometry. A Riemannian metric with this property can be obtained from a $K$-invariant inner product on a complement to $\mathfrak{k}$ in $\mathfrak{u}$ that is compatible with the (almost) complex structure, where $\mathfrak{u}$ and $\mathfrak{k}$ are the Lie algebras of $U$ and $K$, respectively---as we will see, the metrics that occur will in fact satisfy a stronger property---and we assume that we are given such a metric in the following discussion.

As $U$ is a maximal compact subgroup of $G$, given a maximal torus of $U$, its complexification is a maximal torus of $G$, and the respective root systems we obtain are the same.  Furthermore, the set of simple roots yielding $P$ on the algebraic side will then correspond to the compact subgroup $K = P \cap U = L \cap U$.  Now, consider an automorphism $\phi \in E_P$; recall that $\phi$ comes from an automorphism of $G$ as an algebraic group which in turn arises from an automorphism of the Dynkin diagram for $G$ which preserves the nodes corresponding to $P$.  From the construction of the group automorphism induced by a given automorphism of the Dynkin diagram, which is explained for example in \cite[14.2, Theorem]{Humphreys}, from the choice of maximal tori in $U$ and $G$ made here, and from the resulting interplay of root spaces, e.g., see \cite[Ch.\ III, proof of Theorem 6.3]{Helgason}, it follows that $\phi(U) = U$ and $\phi(K) = K$, and hence $E_P \ltimes U$ makes sense as a subgroup of $E_P \ltimes G = \Aut M$.  Since $E_P$ is finite and hence compact, so is $E_P \ltimes U$, and since $U$ is a maximal compact subgroup of $G$, $E_P \ltimes U$ is a maximal compact subgroup of $\Aut M$.  Since $M$ is compact, so is its isometry group \cite[\S5]{MyersSteenrod1939}, which by assumption above, contains $U$ and hence we can conclude the following.

\begin{prop} \label{p:holisoM}
Given a $U$-invariant K\"ahler metric on $G/P = U/K$, we have \begin{align*}
U \leq \holIso(M) \leq E_P \ltimes U.
\end{align*}
\end{prop}

\subsubsection{Projective embeddings} 

Let $M = G/P = U/K$ be a flag variety as above.  Then one has an identification $\Pic M = \mathcal{X}(P)$, the latter being the character group of $P$.  Under this correspondence, the very ample line bundles on $M$ correspond precisely to those characters that are anti-dominant.

Let us fix an anti-dominant character $\chi : P \to L \to \C^\times$; this yields a very ample line bundle $\O_M^\chi(1) = \O_M(1) = G \times^{P, \chi} \C$ and a projective embedding
\begin{align*}
\iota_\chi : M \hookrightarrow \P\left( H^0\big(M, \O_M(1) \big)^\vee \right) = \P^N, \quad \text{ with } N = h^0\big( M, \O_M(1) \big) -1;
\end{align*}
see, e.g., \cite[\S6.1.13]{ChrissGinzburg}.  The $G$-action on $M$ yields a linear one on $H^0\big( M, \O_M(1) \big)$, and hence a morphism of algebraic groups $\theta : G \to SL(N+1)$ fitting into a commutative diagram
\begin{align*}
\xymatrix{ G \ar[r]^-\theta \ar[d] & SL_{N+1} \C \ar[d] \\ M \ar[r]_{\iota_\chi} & \P^N }
\end{align*}
Furthermore, restricting $\theta$ to $U \leq G$, the image $\theta(U)$ will lie in a maximal compact subgroup of $SL_{N+1}(\C)$, which we may conjugate to $SU(N+1)$.  Therefore, we may replace the top row of the above diagram by $U \xrightarrow{\theta} SU(N+1)$.  Since the Fubini--Study metric on $\P^N$ is induced by a $SU(N+1)$-bi-invariant metric on $SU(N+1)$, we can conclude the same about the restricted metric on $M$.

\begin{lem} \label{lem:U-biinvariant_pullback} If a K\"ahler metric on the flag variety $M = U/K$ is obtained by pullback from a Fubini--Study metric via a projective embedding, then it arises from a $U$-bi-invariant metric on $U$.
\end{lem}

\subsection{Lifting actions on tori to complexifications}

Let $T$ be a compact torus of (real) dimension $n$ with the flat metric.  We may write $T = V/\Lambda$, where $V$ is a real vector space of dimension $n$ and $\Lambda \leq V$ is a rank $n$ lattice.  Then, setting $V^\C := V \otimes_\R \C$, $T$ has a natural minimal complexification $T^\C = V^\C / \Lambda$. Of course, $T \hookrightarrow T^\C$ is induced by $V \hookrightarrow V^\C$ and the anti-holomorphic involution is induced by complex conjugation on the $\C$-factor of $V^\C$.  Furthermore, we have a biholomorphism $T^\C \cong (\C^\times)^n$ from the following.

\begin{lem} \label{l:cstar}
Let $\Lambda \subseteq \C^r$ be a rank $r$ lattice.  If $\Lambda \otimes_\Z \C = \C^r$, then $\C^r/\Lambda$ is biholomorphic to $(\C^\times)^r$ under a biholomorphism mapping $\Lambda \otimes_\Z \R$ to $(S^1)^r \subset (\C^\times)^r$.
\end{lem}

\begin{proof}
The assumption that $\Lambda \otimes_\Z \C = \C^r$ means that a $\Z$-basis of $\Lambda$ is also a $\C$-basis of $\C^r$, so by making a linear change of coordinates in $\mathbb{C}^r$, one may assume this is the standard basis. In the corresponding coordinates of $\mathbb{C}^r$, the biholomorphism is induced by the $\Lambda$-invariant holomorphic map
\begin{align*}
(u_1, \ldots, u_n) \mapsto \left( e^{2\pi i u_1}, \ldots, e^{2\pi i u_n} \right). & \qedhere
\end{align*}
\end{proof}

Analogous to Corollary \ref{c:Mactionlift} and Lemma \ref{l:Meqvtsection} above, we have the following. 

\begin{prop} \label{p:torusactionlift}
Let $T$ be as above and suppose $\Gamma$ is a finite group acting freely on $T$ by isometries.  Then there exists a unique holomorphic $\Gamma$-action on $T^\C$ for which the inclusion $T \hookrightarrow T^\C$ is equivariant.  Furthermore, this action is algebraic with respect to the algebraic structure induced by the biholomorphism of Lemma \ref{l:cstar}.
\end{prop}

\begin{proof}
Since $\Gamma$ acts freely on $T$, $X := \Gamma \backslash T$ is a smooth flat Riemannian manifold.  One has $\widetilde{X} = V$ and by the classical Bieberbach theorem, $\pi_1(X) = \Sigma \ltimes M$, where $M \leq V$ is a rank $n$ lattice and $\Sigma \leq O(V)$ is a finite subgroup preserving $M$.  From the sequence of Galois coverings $V \xrightarrow{\Lambda} T \xrightarrow{\Gamma} X$ we have an inclusion of fundamental groups $\Lambda \unlhd \Sigma \ltimes M$ with $\Gamma = (\Sigma \ltimes M)/\Lambda$.  The existence of the $\Gamma$-action on $T^\C$ comes simply from the complexified action of $\Sigma \ltimes M$ on $V^\C$; its uniqueness is a consequence of Lemma \ref{l:uniquelift}. Consider the composition $\Sigma \hookrightarrow \Sigma \ltimes M \to \Gamma$;  its kernel is $\Sigma \cap \Lambda \leq \Sigma \ltimes M$.  As $\Lambda$ is torsion-free and $\Sigma$ is a finite group, this intersection must be trivial, so we may think of $\Sigma$ as a subgroup of $\Gamma$ which acts on $T$.  It therefore preserves $\Lambda$.

To show that the $\Gamma$-action on $T^\C$ is algebraic, we will fix a ($\Z$-)basis $\lambda_1, \ldots, \lambda_n$ of $\Lambda$ to obtain an explicit biholomorphism $T^\C \to (\C^\times)^n$, as in the proof of Lemma \ref{l:cstar} above,
\begin{align*}
\sum_{\alpha=1}^n u_\alpha \lambda_\alpha + \Lambda \mapsto \left( e^{2\pi i u_1}, \ldots, e^{2\pi i u_n} \right),
\end{align*}
and so we can take coordinates $z_\alpha = e^{2\pi i u_\alpha}$ on $T^\C$.

It suffices to show that the $\Sigma$- and $M$-actions on $T$ induced by $\Gamma = (\Sigma \ltimes M)/\Lambda$, are algebraic.  Let $\sigma \in \Sigma$.  Then, since $\sigma$ preserves $\Lambda$ as seen above, with respect to the basis $\lambda = (\lambda_1, \ldots, \lambda_n)$, it has matrix $\sigma^\lambda \in GL_n(\Z)$, say $\sigma^\lambda = (s_{\alpha \beta}^\lambda )$, $s_{\alpha \beta}^\lambda \in \Z$. Then one may verify that $\sigma$ acts on $(\C^\times)^n$ by
\begin{align*}
\sigma \cdot (z_1, \ldots, z_n) = \left( \prod_{\alpha=1}^n z_\alpha^{s_{1 \alpha}^\lambda}, \ldots, \prod_{\alpha=1}^n z_\alpha^{s_{n \alpha}^\lambda} \right).
\end{align*} On the other hand, let $\mu \in M$ and write $\mu = \sum_{\alpha=1}^n m_\alpha \lambda_\alpha$.  Then one may check that the action of $\mu$ on $T^\C$ is given by
\begin{align*}
\mu \cdot (z_1, \ldots, z_n) = \left( e^{2\pi i m_1} z_1, \ldots, e^{2\pi i m_n} z_n \right),
\end{align*}
which is again algebraic in the $z_\alpha$.  In fact, we can say more:  as above, the kernel of $M \hookrightarrow \Sigma \ltimes M \to \Gamma$ is $\Lambda \cap M$, and since the image is finite, one has $[M : \Lambda \cap M] < \infty$.  This implies that $\Lambda \cap M$ is a rank $n$ lattice and hence also $[\Lambda : \Lambda \cap M] < \infty$.  It follows that there exists a $d \in \Z_{> 0}$ such that $d \lambda_\alpha \in M$ for all $\alpha$. Consequently, the $m_\alpha$ above are rational numbers, and the factors $e^{2\pi i m_\alpha}$ hence are, in fact, powers of roots of unity.
\end{proof}

\subsection{Curvature on K\"ahler manifolds}

We follow the conventions of \cite{GoldbergKobayashi1967}.  Let $X$ be a K\"ahler manifold and let $R$ be the Riemannian curvature tensor associated to the K\"ahler metric.  Then we recall that for a (real) $2$-plane $\pi \subseteq T_x X$ in the tangent space to $X$ at a point $x \in X$, the \textbf{sectional curvature of $\pi$} is given by
\begin{align*}
K(\pi) := R(u, v, u, v),
\end{align*}
where $u, v$ form an orthonormal basis of $\pi$.  Restricting this to complex lines in $T_x X$, i.e., subspaces invariant under the complex structure $J$, we get the \textbf{holomorphic sectional curvature}:  if $\sigma$ is such a plane and $u \in \sigma$ a unit vector, then
\begin{align*}
H(\sigma) := R(u, Ju, u, Ju).
\end{align*}
Generalising this slightly is the \textbf{holomorphic bisectional curvature}:  if $\sigma, \sigma'$ are $J$-invariant planes, then this is defined as
\begin{align*}
H(\sigma, \sigma') := R(u, Ju, v, Jv),
\end{align*}
where $u \in \sigma, v \in \sigma'$ are any unit vectors.  Since there is no chance of ambiguity, we will sometimes write ``bisectional curvature'' rather than ``holomorphic bisectional curvature''. 

\begin{rmk} \label{r:secbis}
It is clear to see that holomorphic bisectional curvature determines the holomorphic sectional curvature, preserving signs.  \cite[Equation (4)]{GoldbergKobayashi1967}, which follows from the Bianchi identity, states that
\begin{align*}
R(u, Ju, v, Jv) = R(u, v, u, v) + R(u, Jv, u, Jv)
\end{align*}
so that the sectional curvature also determines the bisectional curvature, also preserving signs.  In particular, if $X$ has non-negative sectional curvature, then it also has non-negative bisectional curvature. With a little bit more work, one can show that in fact holomorphic sectional curvature determines the curvature tensor completely, but we will not need this.
\end{rmk}

Of course, greatly facilitating our main results is Mok's solution to the generalised Frankel conjecture.

\begin{thm}[\protect{\cite[Main Theorem, p.179]{Mok1988}}] \label{t:Mok} Let $(X,h)$ be an $n$-dimensional compact K\"ahler manifold of non-negative holomorphic bisectional curvature and let $(\widetilde{X}, \tilde{h})$ be its universal covering space, equipped with the pullback K\"ahler metric. Then there exist non-negative integers $k$, $N_1, \ldots, N_\ell$ and irreducible compact Hermitian symmetric spaces $M_1, \ldots, M_p$ of rank $\geq 2$, such that as a K\"ahler manifold $(\widetilde{X}, \tilde{h})$ is isomorphic (i.e., isometrically biholomorphic) to
\begin{align} \label{e:Mokdecomp}
(\C^k, g_0) \times (\P^{N_1}, \theta_1) \times \cdots \times (\P^{N_\ell}, \theta_\ell) \times (M_1, g_1) \times \cdots \times (M_p, g_p) ,
\end{align}
where $g_0$ denotes the Euclidean metric on $\C^k$, $g_1, \ldots, g_p$ are canonical metrics on $M_1, \ldots, M_p$, and $\theta_i$, $1 \leq i \leq \ell$, is a K\"ahler metric on $\P^{N_i}$ having non-negative holomorphic bisectional curvature.
\end{thm}

\begin{rmk}[Uniqueness of K\"ahler forms on irreducible compact Hermitian symmetric spaces] \label{r:irredhss}
By the classification of irreducible compact Hermitian symmetric spaces \cite[\S X.6.3]{Helgason}, any such $M$ is a flag variety $G/P$, where $G$ is a simple algebraic group and $P \leq G$ is a maximal parabolic subgroup.  This implies that $H^2(M, \C) = H^{1,1}(M)$ is $1$-dimensional. One may also express $M = U/K$ as a quotient of compact Lie groups, as in Section \ref{s:auts} above; the K\"ahler metric on $M$ in this realisation arises from an $\Ad(K)$-invariant inner product\footnote{In fact, the inner product is induced by the Killing form of $U$.} on $\mathfrak{u}/\mathfrak{k}$.  By the irreducibility assumption there is only one such up to a positive scalar multiple.  On the other hand, any $G$-equivariant holomorphic embedding of $M$ into the complex projective space $\mathbb{P} V$ associated with a unitary $U$-representation space $V$ provides it with a $U$-invariant K\"ahler metric that is induced by a $U$-bi-invariant metric on $U$, as we have seen in Lemma \ref{lem:U-biinvariant_pullback}. This metric will differ from the given one only by a positive scalar multiple.
\end{rmk}

\subsection{Adapted complex structures}

Given a complete real analytic Riemannian manifold $(M, g)$, there exists a unique complex structure on an open neighbourhood of the zero section of the tangent bundle $TM$ such that for a geodesic $\gamma : \R \to M$, the map $\psi_\gamma : \C \to TM$ given by
\begin{align*}
\psi_\gamma( \sigma + i \tau) = \tau \dot{\gamma}(\sigma)
\end{align*}
is holomorphic, where the right side means scalar multiplication.  This complex structure was introduced in \cite{LempertSzoke1991} and \cite{GuilleminStenzel1991}, arising from the study of solutions to the homogeneous Monge--Amp\`ere equation.  There exists an $R \in \R_{> 0} \cup \{ \infty \}$, and a maximal such, such that if
\begin{align*}
T^R M := \{ v \in TM \mid g(v,v) < R^2 \},
\end{align*}
then the complex structure is defined in $T^R M$, see \cite[Theorem 3.1]{LempertSzoke1991}; $T^RM$ is often referred to as the \textbf{Grauert tube of $M$}.  Furthermore, a bound for the sectional curvature of $g$ can be given in terms of $R$ \cite[Theorem 2.4, 4.2]{LempertSzoke1991}; in particular, if $R = \infty$, in which case one says that $M$ has \textbf{entire Grauert tube}, the sectional curvature of $(M, g)$ has to be non-negative.

In the above, we have discussed the adapted complex structure associated to $(M, g)$ as living on the tangent bundle $TM$, but as the metric $g$ always induces a diffeomorphism (indeed, in our case, a real-analytic isomorphism, as $g$ is assumed to be so) from $TM$ to  $T^*M$, we may also think of it as living on the cotangent bundle.  Furthermore, as described in \cite[\S5]{GuilleminStenzel1991}, the natural symplectic form on $T^*M$ is compatible with this complex structure on $T^*M$ in some neighbourhood of the zero section.  We mention this, since in what follows, it will be more natural for us to consider cotangent bundles rather than tangent bundles, and so when we speak of the adapted complex structure on such, this is what we will mean.

\subsection{Big vector bundles}
Modulo our convention for projective bundles, the following is standard.

\begin{defn} \label{d:big}
Let $X$ be a compact complex manifold of dimension $n$.  We say that a line bundle $L$ over $X$ is \textbf{big} if its Kodaira--Iitaka dimension $\kappa(X, L)$ is equal to $n$, see \cite[Definition 2.1.3]{Lazarsfeld}; equivalently, $L$ is big if there exists a $C > 0$ such that $h^0(X, L^{\otimes m}) \geq C \cdot m^n$ for all sufficiently large and divisible $m$, \cite[Lemma 2.2.3]{Lazarsfeld}. 

We say that a vector bundle $E$ over $X$ is \textbf{big} if the line bundle $\O_{\P E^\vee}(1)$\footnote{Recall our convention in Remark \ref{r:Pconv}.} over the projectivisation $\P E^\vee$  is big, see \cite[Example 6.1.23]{Lazarsfeld}.  Note that this is sometimes referred to in the literature as ``L-big'' (e.g., \cite[Definition 1.4]{Jabbusch2009}).  
\end{defn}

\begin{lem} \label{l:etbig}
Let $X$ and $Y$ be compact complex manifolds, $f : Y \to X$ a finite \'etale cover and suppose that $\Theta_X$ is big.  Then $\Theta_Y$ is also big.
\end{lem}

\begin{proof}
The fibre product $\P \Omega_X^1 \times_X Y$ may be identified with $\P f^*\Omega_X^1$, and since $f$ is \'etale, we furthermore have $f^*\Omega_X^1 = \Omega_Y^1$.  So we obtain a cartesian diagram
\begin{align*}
\xymatrix{ \P \Omega_Y^1 \ar[r]^F \ar[d] & \P \Omega_X^1 \ar[d] \\ Y \ar[r]_f & X. }
\end{align*}
In addition, we observe that $F^* \O_{\P \Omega_X^1}(1) = \O_{\P \Omega_Y^1}(1)$.  Consequently, $\O_{\P \Omega_Y^1}(1)$ has a linear subsystem of maximal Kodaira--Iitaka dimension and is therefore big.
\end{proof}

\section{\textbf{Canonical extensions of K\"ahler manifolds}} \label{s:cancxfn}

\subsection{Definition of the canonical extension} \label{ss:cancxfn}

Let $X$ be a K\"ahler manifold with K\"ahler form $\omega$.  Then the corresponding K\"ahler class $[\omega] \in H^{1,1}(X) \cong H^1(X, \Omega_X^1) \cong \Ext^1(\O_X, \Omega_X^1)$ defines an extension of vector bundles \begin{align} \label{e:Omegaext}
0 \longrightarrow \Omega_X^1 \longrightarrow W \longrightarrow \O_X \longrightarrow 0
\end{align}
as a special case of the exact sequence written in Equation~\eqref{e:Vext}.  We will define the \textbf{canonical extension of $X$ with respect to $[\omega]$} as $Z_{[\omega]}$, using the notation of Section \ref{s:affbdl}. In case the K\"ahler class is fixed, the extension is sometimes denoted by $Z_X$.

Of course, the K\"ahler form itself gives a Dolbeault representative for the extension class $[\omega] \in H^1(X, \Omega_X^1)$.  Then upon choosing local primitives $\rho_i \in \A^{1,0}(X_i)$ with $\delbar \rho_i = \omega|_{X_i}$, Lemma \ref{l:affsn}\eqref{l:cansn} implies that these patch together to give a differentiable section 
\begin{align} \label{e:cancxfnsn}
\rho_\omega : X \to Z_{[\omega]}
\end{align}
of $\pi : Z_{[\omega]} \to X$.  The same statement allows us to realise this (in the Dolbeault realisation of the extension given in Equation~\eqref{e:Omegaext}) as the section $(0,1)$ of $\Omega_X^1 \oplus \O_X$, but viewed with the holomorphic structure arising from $\omega$, cf.~ Equation~\eqref{e:Dolstr}).

\begin{thm}[Totally real embedding defined by K\"ahler form] \label{thm:totallyreal} The section $\rho:= \rho_\omega : X \to Z_{[\omega]}$ is a totally real embedding. Moreover, if $\omega$ is real analytic, then $\rho_{\omega}$ is also real analytic. 
\end{thm}

\begin{proof}
This is of course a local question, so we may work over a coordinate chart $U$.  Given a point $p \in U$, as $\omega$ is a K\"ahler form, we may consider holomorphic geodesic coordinates $z_1, \ldots, z_n$  for $\omega$ centred at $p$, e.g., see \cite[Prop.~3.14]{Voisin}, so that in these coordinates $\omega$ has an expression 
\begin{align*} 
\omega = \frac{i}{2} \sum_{p=1}^n dz_p \wedge d\bar{z}_p + O(|z|^2).
\end{align*}
Then an expression for $\rho$ in the same coordinates reads 
\begin{align*}
\rho = \frac{1}{2i} \sum_{p=1}^n \left( \bar{z}_p + t_p(z) + O(|z|^3) \right) dz_p,
\end{align*}
where the $t_p$ are locally defined holomorphic functions.  We think of this as a map $\rho : U \to U \times \C^n$:
\begin{align*}
z \mapsto \left( z, \bar{z}_1 + t_1(z) + O(|z|^3), \ldots, \bar{z}_n + t_n(z) + O(|z|^3) \right).
\end{align*}
With this, it is not hard to compute directly that $J_{U \times \C^n} \circ d\rho_x(T_x U) \cap d\rho_x(T_x U) = \{0\}$; i.e., the image of $\rho_\omega$ is totally real, as claimed.

For the statement about analyticity, one should recall that the proof of the Dolbeault lemma for $(1,1)$-forms involves integrating locally, and if one starts with a real analytic form, integrating will also yield one.  Furthermore, any two primitives differ by an element of the kernel of $\delbar$, hence a holomorphic form.  Thus, \emph{any} primitive will be real analytic if $\omega$ is.
\end{proof}

\begin{rmk} \label{rmk:scalar_multiples_same_bundle}
 If two K\"ahler forms differ by a positive constant, the two resulting canonical extensions are canonically isomorphic as affine bundles, via an isomorphism mapping the image of the canonical section to the image of the canonical section.  If the K\"ahler classes of two K\"ahler forms differ by a positive constant, the resulting canonical extensions are isomorphic as affine bundles, while the K\"ahler forms will, in general, give different totally real sections.
\end{rmk}

\begin{cor}[Minimal complexification defined by K\"ahler form]
If $\omega$ is real analytic, there exists a basis of Stein neighbourhoods $\Omega_\varepsilon \subset Z_{[\omega]}$ of $\rho_\omega(X)$ in $Z_{[\omega]}$ such that each $\Omega_\varepsilon$ is a minimal complexification of $X$ admitting a surjective, holomorphic map $\pi|_{\Omega_\varepsilon}: \Omega_\varepsilon \to X$.
\end{cor}
\begin{proof}
 By Theorem~\ref{thm:totallyreal}, the assumptions of Theorem~\ref{thm:GrauertNirenberg} are fulfilled for $\rho_\omega: X \to Z_{[\omega]}$.
\end{proof}

\begin{rmk}[Real analytic representatives exist] We emphasise that by Proposition~\ref{t:anrep}, every K\"ahler class $a \in H^{1,1}(X)$ contains a real analytic representative $\omega$, to which the above construction can be applied to yield a corresponding minimal complexification.
\end{rmk}

To our knowledge, the above construction  first appeared in somewhat different form in \cite[\S8]{Semmes} in a study of the complex Monge--Amp\`ere equation; this perspective is taken up again in \cite[\S6]{Lempert2017}, which looks at the geometry of spaces of K\"ahler metrics.\footnote{We thank L\'aszl\'o Lempert for bringing this historical point to our attention.}  The spaces $Z_{[\omega]}$ were also considered in \cite[\S2]{Donaldson2002}, where it was shown that $Z_{[\omega]}$ admits a global holomorphic symplectic form $\Theta$ and that the embedding $\rho_{\omega}$ is Lagrangian with respect to the real part of $\Theta$ and symplectic with respect to the imaginary part.

\subsection{Basic properties}
We will see later in Section \ref{s:nonnegbis} that for compact K\"ahler manifolds, non-negative bisectional curvature is a sufficient condition for the canonical extension itself to be a Stein manifold.  It is not true that it is always Stein (see Example \ref{ex:RS}), however, it always possesses the following property of Stein manifolds.

\begin{prop} \label{p:0dim}
Suppose that $(X, \omega)$ is compact K\"ahler and let $Z=Z_{[\omega]}$ be its canonical extension with respect to $[\omega]$.  Then the only compact analytic subsets of $Z$ are $0$-dimensional.
\end{prop}

\begin{proof}
We give the proof assuming that $T \subseteq Z$ is a connected, compact, complex submanifold, the case of an (irreducible) singular analytic subset $T \subseteq Z$ can be handled using a(n embedded) resolution $\widetilde T \to T$ such that $\widetilde T$ is a compact K\"ahler manifold.

Since $Z$ is an affine bundle over $X$, the intersection of each fibre of $\pi : Z \to X$ with $T$ is a compact analytic subspace of an affine space, so is a finite set of points.  Thus the induced proper map $f : T \to X$ is a finite morphism, and therefore an \'etale covering onto its image away from the ramification locus.  Therefore, the pullback $f^* \omega$ is positive away from this proper analytic subset of $T$ and hence, if $k = \dim_\C T$, one has
\begin{equation} \label{eq:integral_positive}
\int_T f^* \omega^k > 0.
\end{equation}

Recall from Lemma \ref{l:extlem}\eqref{l:PWminusPV} that $Z$ is an open complex submanifold of $\P W$. As $X$ is K\"ahler, the latter is a K\"ahler manifold, hence $Z$ itself is K\"ahler, and finally we conclude that $T$ is a compact K\"ahler manifold.  Now, the extension class $[\pi^* \omega] \in H^1(Z, \pi^* \Omega_X^1)$ is trivial, and hence so is the extension class obtained by restriction to $T$; i.e., we have \[ 0= [\pi^*\omega]|_T = [f^* \omega] \in H^1(T, f^* \Omega_X^1).\] Using the map $H^1(T, f^* \Omega_X^1) \to H^{1,1}(T)$ induced by the pullback map $df: f^* \Omega_X^1 \to  \Omega_T^1$ and the Dolbeaut isomorphism, we conclude that $f^* \omega$ is a $\dbar$-exact $(1,1)$-form on $T$. Since $T$ is compact K\"ahler, the $\partial \dbar$-lemma hence implies that $f^* \omega$ is $d$-exact, from which we conclude that $\int_T f^* \omega^k = 0$,
unless $k=0$. From this and Equation~\eqref{eq:integral_positive} above we conclude that $k=0$, as claimed. \end{proof}

The canonical extensions constructed above enjoy the following universal property.
\begin{prop} 
Let $X$ be a K\"ahler manifold with K\"ahler form $\omega$ and let $Z_X = Z_{[\omega]}$ be the corresponding canonical extension. Suppose $T$ is a Stein manifold with a holomorphic map $p : T \to X$.  Then the pullback of the extension~\eqref{e:Omegaext} to $T$ splits, and every splitting induces a holomorphic map $\phi : T \to Z_X$ over $X$, i.e., a holomorphic map making the following diagram commutative,
\begin{align*}
\xymatrix{ 
T \ar[r]^\phi \ar[dr]_p & Z_X \ar[d]^{\pi} \\
& X. }
\end{align*}
\end{prop}

\begin{proof}
We note that if $T$ is Stein, then clearly $H^1(T, f^* \Omega_X^1) = 0$, so that the pullback of the extension class is trivial. The claim therefore is just Lemma \ref{l:extlem}\eqref{l:extlemc} in the case $V = \Omega_X^1$. 
\end{proof}

The following will be useful in determining canonical extensions in particular examples.

\begin{lem} \label{l:cxfnbas}
\begin{enumerate}[(a)] 
\item \label{l:cancxfnprod} Let $(X, \omega_X)$ and $(Y, \omega_Y)$ be K\"ahler manifolds.  Then there is an isomorphism $\alpha : Z_{X \times Y} \arsim Z_X \times Z_Y$ of affine bundles over $X \times Y$ such that the diagram
\begin{align*}
\xymatrixrowsep{1.75pc}\xymatrixcolsep{1pc}
\xymatrix{
Z_{X \times Y} \ar[rr]^\alpha \ar[dr] & & Z_X \times Z_Y \ar[dl] \\
& X \times Y \ar@/^/[ul]^\rho \ar@/_/[ur]_{\rho_X \times \rho_Y} & }
\end{align*}
commutes, where $\rho$ is the section arising from the form $p_X^* \omega_X + p_Y^* \omega_Y$, with $p_X$ and $p_Y$ the respective projections.

\item \label{l:cxfnpullback} Suppose that $\pi : Y \to X$ is a locally isometric holomorphic covering map of K\"ahler manifolds.  Then $Z_Y$ and $Z_X$ are related by a cartesian diagram as follows,
\begin{align*}
\commsq{ Z_Y }{ Z_X }{ Y }{ X }{ \Pi }{}{}{ \pi}.
\end{align*}
In particular, $\Pi$ is also a holomorphic covering.  If $\pi$ is a Galois covering with group $\Gamma$, then so is $\Pi$.  Moreover, the section induced by the K\"ahler form on $Y$ is given by pullback of the section induced by the K\"ahler form on $X$.
\end{enumerate}
\end{lem}

\begin{proof}
The first statement is clear from the fact that the K\"ahler form on $X \times Y$ is the sum of the respective pullback from $X$ and $Y$.  The second comes from Lemma \ref{l:extlem}\eqref{l:baspropb}: since $\pi$ is a covering, $\pi^* \Omega_X^1 = \Omega_Y^1$ and the fact that it is locally isometric means that $\pi^* \omega_X = \omega_Y$.  This equation also yields the last claim of the lemma.
\end{proof}
\subsection{Examples I:  Flat metrics on Euclidean space and complex tori} \label{s:Euctori}

\subsubsection{Euclidean space}

We consider $\C^n$ with the flat metric, so that the associated  K\"ahler form is the standard one given by
\begin{equation} \label{eq:standard_metric}
\omega = \frac{i}{2} \sum_{p=1}^n dz_p \wedge d\bar{z}_p,
\end{equation}
for a choice of global coordinates $z_1, \ldots, z_n$.

We wish to give here a description of the canonical extension $Z = Z_{\C^n}$ of $\C^n$.  Of course, since $H^1(\C^n, \Omega_{\C^n}^1) = 0$, by Lemma \ref{l:affsn}, we will have $Z \cong |\Omega_{\C^n}^1| \cong \C^{2n}$.  However, the isomorphism is non-trivial and will be important when we look at complex tori, so we go through the computation here.  The extension \eqref{e:Vext} for $\Omega_{\C^n}^1$ is $0 \to \Omega_{\C^n}^1 \to W \to \O_{\C^n} \to 0$, 
and over $\C^n$, we have a global frame $dz_1, \ldots, dz_n$ for $\Omega_{\C^n}^1$ and $1$ for $\O_{\C^n}$.  The holomorphic structure on $W$ is given by \eqref{e:Dolstr}:
\begin{align*}
\dbar_W = \mat{ \dbar_{\Omega} }{ \omega }{0 }{ \dbar_{\O}}.
\end{align*}
Therefore, $W$ has a holomorphic frame given by
\begin{align*}
 (dz_p, 0) \quad 1 \leq p \leq n, \quad \left( \frac{i}{2} \sum_{p=1}^n \bar{z}_p \, dz_p, 1 \right).
\end{align*}
Thus, global coordinates on the total space of $W$ are given by
\begin{align*}
(z, w_1, \ldots, w_n, y) \mapsto \left( z, \sum_{p=1}^n \left( w_p + \frac{i}{2} y \bar{z}_p \right) dz_p, y \right).
\end{align*}
The affine bundle $Z = Z_{\C^n}$ is the submanifold where $y = 1$, so this has global coordinates
\begin{align*}
(z,w) \mapsto \left( z, \sum_{p=1}^n \left( w_p + \frac{i}{2} \bar{z}_p \right) dz_p \right).
\end{align*}
If we started with coordinates $(z, u)$ on $| \Omega_{\C^n}^1|$, say
\begin{align*}
(z,u) \leftrightarrow \left( z, \sum_{p=1}^n u_p \, dz_p \right),
\end{align*}
then we have 
\begin{align} \label{e:Cnholcoords}
w_p = u_p - \tfrac{i}{2} \bar{z}_p.
\end{align}
The section of $Z \to \C^n$ given by Lemma \ref{e:cancxfnsn} corresponding to $\omega$ is
\begin{align*}
z \mapsto -\frac{i}{2} \sum_{p=1}^n \bar{z}_p \, dz_p.
\end{align*}
In the holomorphic coordinates $(z,w)$, this is $z \mapsto (z, w = -\tfrac{i}{2} \bar{z})$, and so it is easy to see that this is a totally real embedding.

\subsubsection{Complex tori}

Suppose now that $X$ is a complex torus of dimension $n$ so that $X = \C^n / \Lambda$ for some full rank lattice $\Lambda \cong \Z^{2n}$ in $\mathbb{C}^n$.  We consider $X$ as a K\"ahler manifold by endowing it with flat metric induced from the $\Lambda$-invariant metric \eqref{eq:standard_metric} on $\C^n$.  

\begin{prop} \label{p:Ztorus}
For a complex torus $X$ of dimension $n$ and the flat metric $\omega$, the canonical extension $Z_X = Z_{\omega}$ is biholomorphic to $(\C^\times)^{2n}$.  In particular, it is a Stein manifold.  Moreover, the biholomorphism maps the image of the section induced by the K\"ahler form to the set of real points $(S^1)^{2n} \subset (\mathbb{C}^\times)^{2n}$.
\end{prop}

\begin{proof}
By Lemma~\ref{l:cxfnbas}\eqref{l:cxfnpullback}, the canonical extension $Z_X$ is the quotient of the canonical extension $Z_{\mathbb{C}^n}$ by an induced action of $\Lambda$ which is described as follows.  As the action of $\Lambda$ on $\C^n$ is by translation, the induced action on the cotangent bundle $\Omega_{\C^n}^1$ is simply $\lambda \cdot (z, u) = (z + \lambda, u)$. In terms of the holomorphic coordinates \eqref{e:Cnholcoords} on $Z_{\C^n}$, we get $\lambda \cdot (z,w) = \left(z + \lambda, w - \tfrac{i}{2} \bar{\lambda} \right)$. Hence, this action of $\Lambda$ generates the lattice in $Z \cong \C^{2n}$ given by 
\begin{equation} \label{eq:lifted_lattice}
 \widehat{\Lambda} := \left\{ \left( \lambda, -\tfrac{i}{2} \bar{\lambda} \right) \mid \lambda \in \Lambda \right\}.
\end{equation}

We claim that
\begin{align} \label{e:toruslatticeclaim}
\widehat{\Lambda} \otimes_\Z \C = \C^{2n}.
\end{align}
In order to see this, choose a $\Z$-basis $\lambda_1, \ldots, \lambda_{2n}$ of $\Lambda$ and let $\Pi \in M_{n \times 2n}(\C)$ denote its matrix with respect to the standard basis of $\C^n$.  We note that complex conjugation is with respect to the standard basis in the sense that the standard basis vectors are invariant under conjugation.  By \cite[Proposition 1.1.2]{BirkenhakeLange}, since $\Pi$ arises from a lattice yielding a complex torus, we have
\begin{align*}
\left[ \begin{array}{c} \Pi \\ \overline{\Pi} \end{array} \right] \in GL_{2n}(\C);
\end{align*}
the argument simply uses the fact that the columns of $\Pi$ are linearly independent over $\R$.  But the matrix of $\widehat{\Lambda}$ with respect to the standard basis of $\C^{2n}$ is
\begin{align*}
\left[ \begin{array}{c} \Pi \\ -\frac{i}{2} \overline{\Pi} \end{array} \right] = \mat{ I_n }{ 0 }{ 0 }{ -\frac{i}{2} I_n } \left[ \begin{array}{c} \Pi \\ \overline{\Pi} \end{array} \right] \in GL_{2n}(\C).
\end{align*}
This proves our claim \eqref{e:toruslatticeclaim} and so we may conclude  that $Z_X$ is biholomorphic to $(\C^\times)^{2n}$ by applying Lemma \ref{l:cstar}. 

The set of real points $(S^1)^{2n} \subset (\C^\times)^{2n}$ is precisely the image of $\widehat{\Lambda} \otimes_\Z \R$ in $\widehat{\Lambda} \otimes_\Z \C / \widehat{\Lambda}$.  But then comparing the expression for the section \eqref{e:Cnholcoords} and the definition of $\widehat{\Lambda}$ \eqref{eq:lifted_lattice}, we see that these will have the same image.
\end{proof}

\begin{rmk}
One observes that the biholomorphism $Z_X \cong (\C^\times)^{2n}$ is not algebraic.
\end{rmk}

\subsubsection{Some remarks on abelian varieties}

The results of \cite{BiswasHurtubiseRaina2011} allows us to make the following observation.

\begin{prop} \label{p:abvar}
Let $A$ be a complex abelian variety endowed with the \emph{flat} metric.  Then, as an algebraic affine bundle, $Z_A$ is isomorphic to the de Rham moduli space $\MM_{\dR} = \MM_{\dR}(A, 1)$ of flat rank $1$ connections on the dual abelian variety $A^\vee = \Pic^0 A$.
\end{prop}

\begin{rmk} \label{r:notalgbutStein}
By \cite[Proposition 4.1]{BiswasHurtubiseRaina2011}, $\MM_{\dR}$ admits no non-constant algebraic functions.  This implies that $Z_A$ cannot be an affine variety.  On the other hand, supposing $\dim A = n$, we know that $\MM_{\dR}$ is biholomorphic to the space $\MM_{\B}$, which parametrises $1$-dimensional representations of $\pi_1(A) \cong \Z^{2n}$; thus $\MM_\B \cong (\C^\times)^{2n}$.  Of course, we also knew that $Z_A$ and $(\C^\times)^{2n}$ were biholomorphic from Proposition \ref{p:Ztorus}. In the case $A$ is an elliptic curve, this is the famous example, attributed to Serre, of an algebraic variety which is not affine, but whose underlying complex manifold is Stein.\footnote{This appears to be first written down in detail in \cite[\S6.3, pp.~232--235]{Hartshorne} and is later clarified in \cite[\S7]{Neeman}.}

One could also consider $A$ with a K\"ahler metric arising from a projective embedding (see Section \ref{s:fvHss} below).  It will again be true that the resulting $Z_A$ cannot be an affine variety; indeed, this follows from Corollary \ref{c:big} and the fact that $\Theta_A$ is trivial, as a trivial vector bundle over a base of positive dimension cannot be big.
\end{rmk}

\begin{proof}[Proof of Proposition \ref{p:abvar}]
We may write $A = V / \Lambda$ for some $n$-dimensional vector space $V$ and a full rank lattice $\Lambda \subseteq V$.  The flat metric on $A$ is induced from that on $V$:  if $\varphi_1, \ldots, \varphi_n \in V^\vee$ is a basis of the dual space, then the hermitian metric $h$ is $h = \sum_{\alpha = 1}^n \varphi_\alpha \otimes \overline{\varphi}_\alpha \in V^\vee \otimes \overline{V}^\vee$. Now, it is well-known (e.g., by \cite[Thm.~1.4.1(b)]{BirkenhakeLange}) that we have canonical identifications $H^{1,1}(A) = H^1(A, \Omega_A^1) = H^1(A, \O_A) \otimes V^\vee = \overline{V}^\vee \otimes V^\vee$.
Since $h : V \to \overline{V}^\vee$ gives an isomorphism, we may use it to identify
\begin{align} \label{e:abA}
H^{1,1}(A) = V \otimes V^\vee = \End V,
\end{align}
and under this isomorphism, $h$ corresponds to $\text{Id}_V$.  Typically, we use the corresponding K\"ahler form $\omega = \frac{i}{2} \sum \varphi_\alpha \wedge \overline{\varphi}_\alpha$, but, up to scalars, this also yields the identity.  The Picard variety $\Pic^0 A$ of $A$ can be explicitly described as the dual complex torus $A^\vee = \overline{V}^\vee/ \Lambda^\vee$ \cite[\S2.4]{BirkenhakeLange}, where 
\begin{align*}
\Lambda^\vee = \{ \varphi \in \overline{V}^\vee \mid \langle \varphi, \lambda \rangle \in \Z \quad \forall \lambda \in \Lambda \}.
\end{align*}

Let us consider the moduli space $\MM$ of pairs $(L, \nabla)$, where $L$ is a line bundle on $A^\vee$ and $\nabla$ is a $\lambda$-connection\footnote{We refer the reader to \cite[\S2]{BiswasHurtubiseRaina2011} for a review of $\lambda$-connections.} on $L$.  In the case $\lambda = 0$, we make the further assumption that $L \in A = \Pic^0 A^\vee$; if $\lambda \neq 0$, this must necessarily hold \cite[\S2]{BiswasHurtubiseRaina2011}.  Forgetting the $\lambda$-connection, the map $(L, \nabla) \mapsto L$ 
yields a morphism $\MM \to A$, and this is in fact a vector bundle which is an extension of $\O_A$ by the trivial bundle with fibre $H^0(A^\vee, \Omega_{A^\vee}^1)$. From \cite[(2.4)]{BiswasHurtubiseRaina2011},\footnote{Note that we have interchanged $A$ and $A^\vee$ in their notation.} we see that there is a short exact sequence
\begin{equation} \label{eq:sequenceI}
0 \to A \times H^0( A^\vee, \Omega_{A^\vee}^1) \to \MM \to \O_A \to 0.
\end{equation}
There is a canonical identification $H^0(A^\vee, \Omega_{A^\vee}^1) = \overline{V}$ and again, using the identification $V = \overline{V}^\vee$ coming from the metric $h$, we may also identify $\overline{V}$ with $V^\vee$, hence $H^0(A^\vee, \Omega_{A^\vee}^1) = \overline{V} = V^\vee = H^0(A, \Omega_A^1)$,
and therefore the sequence given in Equation~\eqref{eq:sequenceI} may be identified with
\begin{align} \label{e:abB}
0 \to \Omega_A^1 \to \MM \to \O_A \to 0.
\end{align}
\cite[Lemma 2.1]{BiswasHurtubiseRaina2011} now states that the extension class for this sequence is $-\mathrm{Id}_V$, when we use the identification in \eqref{e:abA}.  Now, the fibre over $1$ in \eqref{e:abB} is the moduli space $\MM_{\dR} = \MM_{\dR}(A, 1)$ of flat rank $1$ connections on $A$, but then this means that, algebraically, $Z_A \cong \MM_{\dR}$. 
\end{proof}

\subsection{Examples II:  Flag varieties} \label{s:fvHss}

We begin here by explaining how, for a smooth projective variety $X$ embedded in $\P^N$, the extension sequence \eqref{e:Omegaext} is related to the Euler sequence for $\P^N$, and to the Atiyah sequence for the $\C^\times$-bundle coming from the line bundle $\O_{\P^N}(1)|_X$.  After this, we specialise the discussion to the case of flag varieties.

Let $X$ be a smooth complex projective variety.  Then upon fixing a projective embedding $X \hookrightarrow \P^N$, $X$ obtains a K\"ahler structure $\omega$ by restricting the Fubini--Study form $\omega_{FS}$ on $\P^N$: $\omega := \omega_{FS}|_X$.  As usual, we let $\O_X(1) := \O_{\P^N}(1)|_X$. In this setup we have the following relation between the canonical extension induced by $\omega$ and the restricted line bundle $\O_X(1)$. 

\begin{lem} \label{p:AtseqO1}
For a smooth projective variety $X$, with K\"ahler structure arising from a projective embedding, the canonical extension \eqref{e:Omegaext} is the dual of the Atiyah sequence associated to the $\C^\times$-bundle underlying the very ample line bundle $\O_X(1)$.
\end{lem}

\begin{proof}
Of course, $\omega_{FS}$ is the curvature form of the standard metric on $\O_{\mathbb{P}^N}(1)$, which restricts to a metric with curvature form $\omega = \omega_{FS}|_X$ on $\O_X(1)$. The Chern connection of the induced metric on $\O_X(1)$ is of type $(1,0)$ and of course yields a connection on the underlying $\C^\times$-bundle with the same curvature form. By \cite[Prop.~4, p.191]{Atiyah}, the extension class of this curvature form is that of the Atiyah sequence, which justifies the claim.
\end{proof}

Now, let $M = G/P = U/K$ be a flag variety; we revert to the notation of Section \ref{s:fvaut}, so that, in particular, $G$, $P$, $L$, $U$, and $K$ have the meaning given there.  Furthermore, we will let $\g$ and $\p$ denote the Lie algebras of $G$ and $P$, respectively. Let us choose a projective embedding $M \hookrightarrow \P^N$; equivalently, we choose a character $\chi : P \to L \to \C^\times$ strictly anti-dominant for $P$ or a very ample line bundle $\O_M(1)$.  As mentioned above, this choice yields a K\"ahler structure on $M$.

\begin{prop} \label{p:fvcxfn}
Let $M = G/P$ be a flag variety endowed with the K\"ahler structure arising from a choice of projective embedding.  Then the canonical complexification $Z_M$ is biholomorphic over $M$ to the smooth affine variety $G/L$, with the projection map from the latter being the canonical one, $G/L \to G/P = M$.  In particular, $Z_M$ is a Stein manifold.  Furthermore, the biholomorphism  may be taken so that the canonical section $M \to Z_M$ corresponds to the natural inclusion $U/K \to G/L$.
\end{prop}

\begin{rmk}
If $a, a' \in H^1(M, \Omega_M^1)$ are non-proportional ample classes arising from (necessarily different) equivariant projective embeddings of $M$, the corresponding spaces $Z_a, Z_{a'}$ are not isomorphic as affine bundles.  Thus, if we take the composition of the biholomorphism $Z_a \to G/L$ provided by Proposition~\ref{p:fvcxfn} with the inverse of $Z_{a'} \to G/L$, while fibre-preserving, it will not in general preserve the affine-linear structure of the fibres.
\end{rmk}

\begin{proof}[Proof of Proposition~\ref{p:fvcxfn}]
By Proposition \ref{p:AtseqO1} the canonical extension \eqref{e:Omegaext} for $M$ is obtained as the dual of the Atiyah sequence for $Q$.  It is explained in \cite[\S4.3]{BiswasWong2017} that this Atiyah sequence is the sequence of vector bundles associated to a sequence of $P$-representations
\begin{align} \label{e:Prep}
0 \to (\g/\p)^\vee \to \mathfrak{w} \to \C \to 0,
\end{align}
where $\C$ denotes the trivial $1$-dimensional $P$-representation.  If $\a \subseteq \mathfrak{w}$ is the preimage of $1 \in \C$, then $\a$ is a $P$-invariant affine subspace of $\mathfrak{w}$ modelled on the $P$-representation $(\g/\p)^\vee$ and 
\begin{align*}
Z_M = G \times^P \a.
\end{align*}
Furthermore, there is an element $\nu_0 \in \a$ which has stabiliser precisely $L$, and the $G$-equivariant morphism $G/L \to Z_M = G \times^P \a$ given by $g \mapsto [g, \nu_0]$ is an isomorphism, see \cite[Proposition 4.20]{BiswasWong2017}.  Of course, the projection map $Z_M \to M$ is simply the canonical projection $G/L \to G/P$.

Now, if we use the realisation $M = G/P = U/K$ of Section \ref{s:auts}, then considering \eqref{e:Prep} as a sequence of $K$-representations, noting $K \leq P$, the corresponding sequence of vector bundles associated to the $K$-bundle $U \to M = U/K$ yields the ($C^\infty$-version) of \eqref{e:Omegaext}  for $M$.  But now, since $K$ is compact, \eqref{e:Prep} splits as a sequence of $K$-modules.  On the other hand, we obtain a section of $Z_M$ from the (totally real) inclusion 
\begin{equation} \label{e:realsnM}
 U/K \hookrightarrow G/L
\end{equation}
 given by \eqref{eq:section_in_GmodL}. Now, under the $K$-splitting $\mathfrak{w} = (\g/\p)^\vee \oplus \C$, since $\nu_0 \in \a$, the image of this section in $\C$ is $1$, so it must go to an element of the form $(\tau_0, 1) \in (\g/\p)^\vee \oplus \C$ for some $\tau_0 \in (\g/\p)^\vee$.  This $\tau_0$ yields a (real analytic) section of $U \times^K (\g/\p)^\vee \cong \Omega_M^1$, hence $\tau_0$ is a $K$-invariant element of $(\g/\p)^\vee$.  But since $\chi$ was chosen to be strictly anti-dominant, $(\g/\p)^\vee$ has no trivial sub-$K$-representations, hence $\tau_0 = 0$.

It follows that the section $M \to Z_M$ induced from the inclusion \eqref{e:realsnM} corresponds to the section $(0,1)$ of $\Omega_M^1 \oplus \O_M$, under the  splitting arising from the $K$-splitting of \eqref{e:Prep}.  Thus, from Lemma \ref{l:affsn}(\ref{l:cansn}), this yields the canonical section (arising from the K\"ahler form).
\end{proof}

Together with the discussion in Remark \ref{r:irredhss} the previous result yields

\begin{cor}\label{cor:canonical_complexification_IHSS}
Let $M = G/P$ be an irreducible Hermitian symmetric space with the symmetric metric. Then, the canonical complexification $Z_M$ is  biholomorphic  over $M$ to the smooth affine variety $G/L$, with the projection map being the canonical one $G/L \to G/P = M$.  In particular, $Z_M$ is a Stein manifold.  Furthermore, the isomorphism can be chosen so that the canonical section $M \to Z_M$ corresponds the natural inclusion $U/K \hookrightarrow G/L$.
\end{cor}

\section{Non-negative holomorphic bisectional curvature} \label{s:nonnegbis}

In this section, we consider the canonical extensions for compact K\"ahler manifolds of non-negative bisectional curvature. With the help of a structural result, which we prove first, we can relate this class of manifolds and hence their extensions to the ones  already considered in previous sections.

\subsection{Structure theory}
The rough structure of compact K\"ahler manifolds with non-negative bisectional curvature is given by the following result, cf.~\cite{CDP}. 

\begin{prop} \label{p:bisfinitecover}
Let $(X, \omega)$ be a compact K\"ahler manifold of non-negative bisectional curvature.  
\begin{enumerate}[(a)]
\item \label{p:bisfinitecovera} Let $\pi := \pi_1(X)$.  Then if we take the universal cover $\widetilde{X}$ as in \eqref{e:Mokdecomp}, $\pi$ acts freely on the factor of $\C^k$, so that $R := \C^k/\pi$ is a flat (compact) K\"ahler manifold.  If we let $M$ be the product of the other factors, so that $M$ is a flag manifold carrying a K\"ahler metric of non-negative bisectional curvature, then 
\begin{align*}
X = \C^k \times^\pi M
\end{align*}
which is a locally trivial fibre bundle over $R$ with locally constant transition functions acting by K\"ahler isometries on $M$.

\item \label{p:bisfinitecoverT} There is a finite Galois cover of $R$ by a (flat) complex torus $T$, say with $\pi_1(T) =: \Xi$ a lattice in $\C^k$, with $\Xi$ acting on $M$ via the connected component of its K\"ahler isometry group, such that we have a cartesian diagram
\begin{align}\label{eq:structure_diagram}
\begin{gathered}
\xymatrix{
\C^k \times^\Xi M \ar[r] \ar[d] & X \ar[d] \\
T \ar[r] & R.}
\end{gathered}
\end{align}
In particular, $\C^k \times^\Xi M \to X$ is a finite Galois cover (with the same group as $T \to R$).
\end{enumerate}
\end{prop}

\begin{rmk} \label{r:bisfinitecovertriv}
It is furthermore true that the fibre product $T \times_R X = \C^k \times^\Xi M$ is diffeomorphic (indeed, real-analytically isomorphic) to the product $T \times M$.  In fact, in the proof of Theorem \ref{thm:good_complexification} we will construct an explicit such isomorphism (depending only on the action of $\Lambda$ on $M$), so a justification of this fact will be deferred until then.
\end{rmk}

\begin{proof}
In the following, for a K\"ahler manifold $(Y, g_Y)$ we denote by $\holIso(Y)$ the corresponding group of holomorphic isometries (i.e., the intersection of the holomorphic automorphism group with the isometry group).  

Our first observation is that if $g$ is any K\"ahler metric on a simply-connected, compact complex manifold $Y$ with $h^{1,1}(Y) = 1$, then the holonomy representation of $(Y, g)$ is irreducible.  Otherwise, the de Rham decomposition would give a non-trivial factorisation $(Y, g) \cong (Y_1, g_1) \times (Y_2, g_2)$, via a holomorphic isometry, with both $Y_1$ and $Y_2$ compact (see \cite[\S1 Th\'eor\`eme, p.\ 756]{Beauville1983} and following remark).  The latter statement about compactness implies that $h^{1,1}(Y_i) \neq 0$, since the K\"ahler metric $g_i$ yields a non-trivial class in $H^{1,1}(Y_i)$; but then one would have $h^{1,1}(Y) = h^{1,1}(Y_1) + h^{1,1}(Y_2) \geq 2$, contradicting our assumption.  In particular, by Remark \ref{r:irredhss}, each of the non-Euclidean factors in the decomposition \eqref{e:Mokdecomp} provided by Theorem \ref{t:Mok} has irreducible holonomy, and hence \eqref{e:Mokdecomp} is a de Rham decomposition.

Our argument now follows the proof of the Bogomolov--Beauville structure theorem for compact K\"ahler manifolds with vanishing first Chern class \cite[Th\'eor\`eme 1]{Beauville1983}.  Let $M$ denote the compact factor $\prod_i \P^{N_i} \times \prod_j M_j$ in \eqref{e:Mokdecomp}, together with the metric on it described by Theorem \ref{t:Mok}, and consider the action of $\pi := \pi_1(X)$ on $\widetilde{X} = \C^k \times M$.  By the uniqueness of the de Rham decomposition (see \cite[Remarque, p.~757]{Beauville1983}), one concludes that every element of $\holIso(\C^k \times M)$ is of the form $g \times h$ for $g \in \holIso(\C^k)$, $h \in \holIso(M)$.  We may therefore realise $\pi$ as a subgroup of $\holIso(\C^k) \times \holIso(M)$.

Let $\Sigma \leq \pi$ be the subgroup of elements of the form $\text{Id}_{\C^k} \times h$ for some $h \in \holIso(M)$, i.e., $\Sigma$ is the kernel of the action of $\pi$ on $\C^k$ via the homomorphism $\holIso(\C^k) \times \holIso(M) \to \holIso(\C^k)$. Then we note that $\pi/\Sigma \leq \holIso(\mathbb{C}^k)$ acts on $\C^k$ with compact quotient, so by the classical Bieberbach theorem, $\pi/\Sigma \cong \Lambda \rtimes \Phi$, where $\Lambda \leq \C^k$ is a full rank lattice (i.e., $\Lambda \cong \Z^{2k}$ and $\Lambda \otimes_\Z \R \cong \C^k$) and $\Phi$ is a finite group of isometries.  Furthermore, $\Sigma$ must act freely on $M$; but by Proposition \ref{p:Hssbas}, no non-trivial group does so.  It follows that $\Sigma = \{ 1 \}$ and hence $\pi \cong \Lambda \rtimes \Phi$. Therefore, $\pi$ acts freely on $\C^k$, which proves \eqref{p:bisfinitecovera} as well.

Take a basis $\lambda_1, \ldots, \lambda_{2k}$ of $\Lambda$.  These act on $M$ via the compact group $\holIso(M)$; since the quotient $\holIso(M)/ \holIso(M)^\circ$ is a finite group, there are $d_i \in \Z_{>0}$, $1 \leq i \leq 2k$, such that $d_i \lambda_i \in \holIso(M)^\circ$.  Consider the sub-lattice $\Lambda' \leq \Lambda$ spanned by the $d_i \lambda_i$; then taking the intersection of all the $\Phi$-orbits of elements in $\Lambda'$, we obtain a $\Phi$-invariant sub-lattice $\Xi \leq \Lambda' \leq \Lambda$.  By construction, $\Xi$ acts on $M$ via $\holIso(M)^\circ$.  It is straightforward to see that diagram\ \eqref{eq:structure_diagram} is indeed cartesian.  Thus, \eqref{p:bisfinitecoverT} is proved.
\end{proof}

\subsection{Steinness of canonical extensions}

Let $M = G/P = U/K$ be a flag variety endowed with an invariant metric coming from an equivariant projective embedding, let $T$ be a complex torus endowed with a flat metric, let $\Xi = \pi_1(T)$ be its fundamental subgroup, and suppose that $\Xi$ acts on $M$ via a homomorphism to $U = \holIso(M)^\circ \leq G = \Aut(M)$.  Then the associated bundle 
\begin{align*}
Y := \C^k \times^\Xi M
\end{align*}
is a compact K\"ahler manifold which is the total space of a fibre bundle over $T$ with fibre $M$ and whose transition functions may be taken to be locally constant.

Using the structural result obtained in the previous section as well as the information already gathered on canonical extensions of tori and flag manifolds, we are now in a position to prove the following.

\begin{thm} \label{t:bisholStein}
Suppose a compact K\"ahler manifold $X$ admits a finite Galois covering isometrically biholomorphic to a fibre bundle $Y = \C^k \times^\Xi M$ of the kind described in the preceding paragraph.  Then, the canonical extension $Z_X$ is Stein.  In particular, the conclusion holds for a compact K\"ahler manifold $X$ endowed with a K\"ahler metric of non-negative bisectional curvature.
\end{thm}

\begin{proof}
By Lemma \ref{l:cxfnbas}\eqref{l:cxfnpullback}, since the bottom row of diagram
\begin{align*}
\commsq{Z_{\C^k \times M} = Z_{\C^k} \times Z_M}{ Z_Y }{ \C^k \times M }{ Y }{}{}{}{}
\end{align*}
is covering with Galois group $\Xi$, so is the top row.  Since $\Xi$ acts on $Z_{\C^k}$ freely with quotient $Z_T$ by Lemma\ \ref{l:cxfnbas}(\ref{l:cxfnpullback}),  we see that $Z_Y$ is a fibre bundle over $Z_T$ with fibre $Z_M$.  Now, we may consider $Z_{\C^k}$ as a principal $\Xi$-bundle over $Z_T$.  Since the action of $\Xi$ on $M$ is given by a homomorphism $\Xi \to U \leq G$, we may consider the principal $G$-bundle over $Z_T$ obtained by extension of structure group $P := Z_{\C^k} \times^\Xi G$.  Then one has $Z_Y \cong P \times^G Z_M$; hence we have a holomorphic fibre bundle whose structure group $G$ is a connected complex Lie group and a whose base $Z_T$ and fibre $Z_M$ are Stein (Propositions \ref{p:Ztorus} and \ref{p:fvcxfn}).  We can therefore conclude that the total space $Z_Y$ is Stein by applying the classical theorem of Matsushima--Morimoto, \cite[Th\'eor\`eme 6]{MatsushimaMorimoto1960}. 

Again by Lemma \ref{l:cxfnbas}\eqref{l:cxfnpullback}, since $Y \to X$ is a finite Galois covering, so is $Z_Y \to Z_X$.  Since $Z_Y$ is Stein, we can conclude the same of $Z_X$, for example, by averaging an exhaustion function on $Z_Y$ so that it descends to $Z_X$.  This proves the first part of the statement.

With the assumption that $X$ is of non-negative bisectional curvature, Proposition \ref{p:bisfinitecover}\eqref{p:bisfinitecoverT} says that $X$ has a finite Galois cover of the form $Y = \C^k \times^\Xi M$.  However, the metric on $M$ is not necessarily induced from a projective embedding.  From Theorem \ref{t:Mok}, we know that $M = \prod M_i$ is a product of irreducible flag varieties $M_i$, each carrying a metric of non-negative bisectional curvature.  If $M_i$ is an irreducible compact Hermitian symmetric space of rank $\geq 2$,  Remark \ref{r:irredhss} implies that the corresponding canonical extension \eqref{e:Omegaext} is isomorphic to that obtained from a projective embedding as in Section \ref{s:fvHss} above.  Furthermore, if $M_i$ is a projective space $\P^N$ with \emph{some} metric of non-negative  bisectional curvature, then since $h^{1,1}(\P^N) = 1$, the extension classes corresponding to that metric and the Fubini--Study metric yield isomorphic canonical complex extensions by Remark\ \ref{rmk:scalar_multiples_same_bundle}; therefore, we may as well assume that the metric on $M$ comes from a projective embedding, since we are only concerned with the biholomorphism type and the action of $G$ on $Z_M$.  We may then conclude by the arguments above.
\end{proof}

\begin{rmk}
Another possible proof of the second part of the above theorem, where the bisectional curvature of $X$ is assumed to be non-negative, was suggested to us by L\'aszl\'o Lempert.  The metric on $X$ induces one on the contangent bundle and hence defines a length function $h$ on the total space $T^*X$.  One can then compute $\partial \dbar$ of this function with respect to the complex structure yielding $Z$ and one sees that it is strictly plurisubharmonic on a neighbourhood of the section $\rho$ (which is the zero section upon identifying $Z = T^*X$ as smooth manifolds).  Furthermore, $h$ will be strictly plurisubharmonic on the entirety of $Z$ if and only if $X$ has non-negative bisectional curvature.  Thus, one would obtain a strictly plurisubharmonic exhaustion function on $Z$ and could thereby conclude it is Stein. 
\end{rmk}

The following example shows that without the curvature condition the canonical extension is not necessarily Stein.

\begin{ex} \label{ex:RS}
Let $X$ be a compact Riemann surface, $V$ a line bundle on $X$ and $a \in H^1(X, V)$.  Consider the extension \eqref{e:Vext} and let $S := \P W$.  Then the ruled surface $S \to X$ admits a section $\sigma$ with image $X_\infty$ at infinity, corresponding to the inclusion $\P V \subseteq S$ or equivalently, to the surjection $W^\vee \twoheadrightarrow V^\vee$.  Clearly, this section meets each fibre once, and hence the line bundle corresponding to the divisor $X_\infty$ is $\O_S(1)$; moreover, we have $\sigma^*(\O_S(1)) = V^\vee$,  so that the self-intersection $(X_\infty)^2 = - \deg V$; cf.~\cite[\S5, Lem.~10]{Friedman}.  In the case that $g(X) > 1$, $V = \Omega_X$ and $a = [\omega] \in H^1(X, \Omega_X)$ is the class of any K\"ahler metric on $X$ (note that $h^1(X, \Omega_X) = 1$ by Serre duality), the self-intersection of $X_\infty$ is negative by the preceding computations, and therefore by Grauert's criterion \cite[\S3, Thm.~20]{Friedman}, $X_\infty$ can be contracted to a point on a normal, compact complex space $\widetilde{S}$ of dimension two.  But then $Z_X = S \setminus X_\infty = \widetilde{S} \setminus \{ \textnormal{pt} \}$.  But the latter admits no non-constant global holomorphic functions. Hence $Z_X$ cannot be Stein.
\end{ex}

\begin{rmk}
It is an open question, motivated by Jouanolou's trick in algebraic geometry and raised for example by Campana and Winkelmann in \cite[\S3]{CampanaWinkelmann}, whether every compact K\"ahler manifold admits a holomorphic surjection $p: \Omega \to X$ from a Stein manifold $\Omega$ whose fibres are individually biholomorphic to $\mathbb{C}^{N}$ (together with a real analytic section). Theorem \ref{t:bisholStein} provides a very modest extension of the Jouanolou result to some non-algebraic cases.
\end{rmk}

\subsection{Existence of good complexifications} \label{s:goodcxfn}

As mentioned in the Introduction, the question of having affine canonical extensions is related to the existence of good complexifications.  With the help of a non-trivial result of Totaro's \cite[Lemma 3.1]{Totaro2003}, we can show that manifolds of a certain class, which includes compact K\"ahler manifolds of non-negative bisectional curvature, indeed admit such good complexifications. 

Let $B$ be a smooth compact manifold with universal covering space $\R^r$ and fundamental group $\pi := \pi_1(B)$ of the form $\Gamma \ltimes \Lambda$, where $\Lambda \leq \R^r$ is a rank $r$ lattice and $\Gamma \leq GL_\Z(\Lambda) \leq GL_r(\R)$ is a finite subgroup.  Let $M$ be a smooth compact manifold and $K \leq \textnormal{Diff}(M)$ a compact (but not necessarily connected) Lie subgroup of diffeomorphisms.  Let $\rho : \pi \to K$ be a homomorphism and let 
\begin{align} \label{e:Gammaaction}
X & := \R^r \times^\rho M = (\R^r \times M)/ \pi, &\text{where } &  & \sigma \acts (x, m) & = \left( \sigma \acts x, \rho(\sigma)(m) \right)
\end{align}
and $\sigma \in \pi = \Gamma \ltimes \Lambda$ acts as an affine transformation on $\R^r$, which gives meaning to $\sigma \acts x$, when $x \in \R^r$.  In particular, if $\lambda \in \Lambda$, then it acts on $\R^r$ by translation, so we will often write the action as $\lambda \acts (x, m) = \left( x + \lambda, \rho(\lambda) (m) \right)$.  Thus, $X$ is a locally trivial fibre bundle over $B$ with locally constant transition functions with fibre $M$.  Of course, $\Gamma$ itself acts on $M$, being a subgroup of $\pi$.

It follows from Proposition \ref{p:bisfinitecover} and its proof as well as from Proposition \ref{p:holisoM} that if $X$ is a compact K\"ahler manifold of non-negative bisectional curvature, then it is of the form just described, with $M$ a flag manifold.  

\begin{thm} \label{thm:good_complexification} 
Let $X$ be as in \eqref{e:Gammaaction} and suppose that $M$ admits a good complexification $W$ such that the $K$-action on $M$ extends to an algebraic action of $K^\C$ on the complex affine variety $W(\C)$.  Then $X$ admits a good complexification.  In particular, this holds if $X$ is a compact K\"ahler manifold with non-negative bisectional curvature.
\end{thm}

To prove this, we apply the following result of Totaro.

\begin{lem}[\protect{\cite[Lemma 3.1]{Totaro2003}}] \label{l:tot}  Let $U$ be a good complexification with an algebraic action of the complexification $G$ of a compact Lie group $G(\R)$.  If $G(\R)$ acts freely on $U(\R)$, then the algebraic group $G$ acts freely on $U$, and the quotient variety $U/G$ is a good complexification of the quotient manifold $U(\R)/G(\R)$.
\end{lem}

\begin{proof}[Proof of Theorem~\ref{thm:good_complexification}]
To show that a compact K\"ahler manifold of non-negative bisectional curvature is a particular case of the first statement, we first observe that if $M = G/P$ is a flag variety, then $W := G/L$ is a good complexification of $M$ (notation as in Section \ref{subsubsect:flag_automo}); this is already \cite[Theorem 5.1]{Kulkarni1978}.  However, it can also be seen from Lemma \ref{l:tot} above, recalling that a flag variety can always be obtained as a quotient of compact Lie groups (\S\ref{s:auts}), and that every compact Lie group is a real algebraic variety (e.g., see \cite[Theorem 5.2.12]{OnishchikVinberg}).  We also need to see that the $\Gamma$-action on $M$ extends to an algebraic one on $W$, but this becomes clear upon considering Theorem \ref{t:AutM} and Proposition \ref{p:holisoM}.

We now prove the more general statement.  Let $Y := \R^r \times^{\Lambda, \rho} M$, so that $Y$ is a finite Galois cover of $X$ with group $\Gamma$.  By Lemma \ref{l:tot}, it then suffices to show that $Y$ admits a good complexification $U$ on which $\Gamma$ acts algebraically.  

Now, $Y$ is a locally trivial fibre bundle over the torus $T := \R^r/\Lambda$ with fibre $M$.  The
next thing to say is that for such a $T$, $(\C^\times)^r$ is a good complexification, as the latter can be identified with the $\C$-points of the real affine variety 
\begin{align*}
V := \Spec \R[x_1, y_1, \ldots, x_r, y_r]/ ( x_1^2 + y_1^2 - 1, \ldots, x_r^2 + y_r^2 - 1);
\end{align*}
of course, $T$ is diffeomorphic to the manifold of real points.

As a point of notation, we will simply write $W$ to mean the complex affine variety $W(\C)$ as well as the corresponding complex manifold.  Also, if $M \to W$ is written as an inclusion, then we implicitly mean $W(\R) \to W(\C)$.  The same goes for $T$ and $V$.

Our candidate $U$ for a good complexification of $Y$ will be $U := V \times W$.  Therefore it remains to construct a diffeomorphism $Y \xrightarrow{\sim} (V \times W)(\R) = V(\R) \times W(\R) = T \times M$ and show that the $\Gamma$-action extends to an algebraic one on $V \times W$.  This action will not be the diagonal one: we will now construct the diffeomorphism explicitly as the restriction of a biholomorphism of corresponding ``complexified'' objects, and transport the action through this biholomorphism in order to show that it is indeed algebraic.

Viewing $\R^r$ as a principal $\Lambda$-bundle over $T$, $Y$ is then the fibre bundle associated to the homomorphism $\rho : \Lambda \to K \leq \textnormal{Diff}(M)$.  Since $\Lambda$ is abelian, $\rho$ must factor through some torus $J \leq K$, therefore we may form the principal bundle $A := \R^r \times^{\Lambda, \rho} J$, and it is clear that $Y = A \times^J M$.

For the ``complexified'' version, we again have fibre bundles over the base $\C^r / \Lambda$ and we observe that it is straightforward to write down a biholomorphism $\C^r /\Lambda \to V$ so that
\begin{align*}
\xymatrix{ \R^r/\Lambda \ar[r] \ar[d] & T \ar[d] \\ \C^r / \Lambda \ar[r] & V }
\end{align*}
commutes; we will thus often identify $\C^r / \Lambda = V$.  The ``complexified'' version of $Y$ will be $Z := \C^r \times^{\Lambda, \rho} W$.\footnote{Writing ``complexified'' in quotation marks, while it is true that the objects here are complexifications in the sense of Definition \ref{d:cxfn}, we will not make use of this; rather, we are only concerned with their structure as complex manifolds.  However, describing them as such helps to convey the idea of the argument here.}  Observe that since $\pi$, and so $\Gamma$, acts on $\R^r$ via affine transformation \eqref{e:Gammaaction}, there is an obvious extension of the action to $\C^r$.  Furthermore, by assumption $K^\C$ acts algebraically on $W$, extending the $K$-action on $M$, and since $\rho : \pi \to K \leq K^\C$, the $\pi$- and $\Gamma$-actions on $W$ are algebraic.  Thus, we obtain an action of $\Gamma$ on $Z$, given by the same expression as in \eqref{e:Gammaaction}.  Let $S := J^\C$ and set $Q := \C^r \times^{\Lambda, \rho} S$, which is a principal $S$-bundle over $V$ (so that this is the ``complexified'' version of $A$).  Explicitly, $Q$ is the quotient of $\C^r \times S$ by the $\Lambda$-action on $\C^r \times S$ given by
\begin{align*}
\lambda \acts (v, s) = \big(v + \lambda, \rho(\lambda)s \big).
\end{align*}
We denote the class of $(v, s) \in \C^r \times S$ in $Q$ by $[v,s]$.  Using the inclusions $\R^r \hookrightarrow \C^r$, $J \leq S$ and $M \hookrightarrow W$, we have inclusions
\begin{align} \label{e:incs}
A & \hookrightarrow Q & \text{and}& & Y & \hookrightarrow Z,
\end{align}
noting that $\R^r \times J$ is a $\Lambda$-invariant subspace of $\C^r \times S$.  With these relations, we will use the same notation for points of $A$ that we introduced above for $Q$.  As above, since $S = J^\C \leq K^\C$, $S$ acts algebraically on $W$, and one has $Z = Q \times^S W$.  

Given any holomorphic map $f : \C^r \to S$ satisfying
\begin{align} \label{e:hp}
f(v + \lambda) \rho(\lambda) = f(v)
\end{align}
for all $v \in \C^r$, $\lambda \in \Lambda$, we can write down an explicit holomorphic trivialisation $\Phi_f : Q \to V \times S$ as follows
\begin{align}\label{eq:isoI}
[v, s] \mapsto ( [v], f(v) s),
\end{align}
where one uses the property \eqref{e:hp} to show that this is well-defined; the inverse $\Phi_f^{-1}$ is given by
\begin{align}\label{eq:isoII}
([v], s) \mapsto [v, f(v)^{-1} s].
\end{align}
Furthermore, it is easy to see that if $f$ satisfies $f(\R^r) \subseteq J$, then $\Phi_f$ restricts to a real-analytic isomorphism $A \xrightarrow{\sim} T \times J$.  Therefore, given such an $f$, we may construct a commutative diagram
\begin{align*}
\xymatrix{
A \ar[r]^-{\Phi_f} \ar[d] & T \times J \ar[d] \\ Q \ar[r]_-{\Phi_f} & V \times S, }
\end{align*}
where the vertical map on the left is the inclusion on the left side of \eqref{e:incs}, and the one on the right is the inclusion $(V \times S)(\R) \subseteq (V \times S)(\C)$.

To show that an $f$ satisfying \eqref{e:hp} always exists, let us assume for the moment that $S = \C^\times$, so that $J = S^1 \leq \C^\times$.  Let $\lambda_1, \ldots, \lambda_r$ be a basis of $\Lambda \leq \mathbb{R}^r$; since $\rho(\Lambda) \leq J$, we may write $\rho(\lambda_j) = \exp(t_j)$ for some $t_j \in i \R \subseteq \C$ for $1 \leq j \leq r$.  Now, if $v := \sum_{j=1}^r v_k \lambda_k \in \C^r$ for some $v_k \in \C$, the map given by
\begin{align*}
f \left( v \right) := \exp\left( -\sum_{l=1}^r t_k v_k \right)
\end{align*}
satisfies \eqref{e:hp} above.  It is also easy to see that if $v \in \R^r$, i.e., $v_k \in \R$ for $1 \leq k \leq r$, then the right hand side takes values in $S^1 = J$.  Furthermore, we observe that $f|_{\Lambda} = \rho^{-1}$.  

For an arbitrary complex algebraic torus $S$, upon choosing an isomorphism $S \cong (\C^\times)^p$, we may repeat the above argument component-wise, and thus obtain an $f$ satisfying \eqref{e:hp} for any $S$ and any $\rho$ with $f(\R^r) \leq J$.  Additionally, the $f$ constructed in this way has the property that it is a homomorphism of complex Lie groups $\C^r \to S$ extending $\rho^{-1}$ on $\Lambda$.  This property will be used in a crucial way below.

Now, $\Phi_f$ induces an isomorphism $\phi_f : Z \to V \times M$, defined as the composition
\begin{align*}
Z = \C^r \times^{\Lambda, \rho} W \xrightarrow{\sim} Q \times^S W \xrightarrow{\Phi_f \times^S \text{Id}_W} (V \times S) \times^S W \xrightarrow{\sim} V \times W;
\end{align*}
explicitly, $\phi_f$ and $\phi_f^{-1}$ are respectively given by
\begin{align*}
[v,w] & \mapsto \left( [v], f(v) (w) \right) & \text{and}& & ([v], w) & \mapsto \left[ v, \left( f(v)^{-1} \right) (w) \right],
\end{align*}
cf.\ \eqref{eq:isoI} and \eqref{eq:isoII}.  Furthermore, using the property that $\Phi_f$ restricts to an isomorphism $A \to T \times J$, it is easy to see that $\phi_f$ restricts to a real-analytic isomorphism $Y \to T \times M$ making the diagram
\begin{align*}
\xymatrix{ Y = \R^r \times^{\Lambda} M \ar[r]^-{\phi_f} \ar[d] & T \times M \ar[d]^\iota \\
Z = \C^r \times^\Lambda W \ar[r]_-{\phi_f} & V \times W }
\end{align*}
commute.\footnote{This proves the statement made in Remark \ref{r:bisfinitecovertriv}.}  In particular, we see that $V \times W$ is a good complexification of $Y$, as $\phi_f$ gives a diffeomorphism (indeed, real-analytic isomorphism) onto the submanifold of real points.

What remains is to show that the $\Gamma$-action on $V \times W$ is algebraic.  This action arises from the action on $Z = \mathbb{C}^r \times^{\Lambda, \rho} W$, which, as we have mentioned, is given by \eqref{e:Gammaaction}.  Given $\gamma \in \Gamma$, we can transport its action to $V \times W$ via $\phi_f$, so that 
\begin{align*}
\gamma \acts \left( [v], w \right) & = \phi_f \left( \gamma \acts \phi_f^{-1}( [v], w) \right) = \phi_f \left( \gamma \acts \left[ v, \left( f(v)^{-1} \right) (w) \right] \right) = \phi_f \left( \left[ \gamma \acts v, \left( \rho(\gamma) f(v)^{-1} \right) (w) \right] \right) \\
& = \left( [ \gamma \acts v] , \left( f( \gamma \acts v) \rho(\gamma) f(v)^{-1} \right) (w) \right) = \left( [ \gamma \acts v] , \left( f( \gamma \acts v) \rho(\gamma) f(v)^{-1} \rho(\gamma)^{-1} \right) \circ \rho(\gamma) (w) \right)
\end{align*} 

The action on the $V$-factor is obtained by descending a linear map on $\C^r$, so is algebraic by Proposition \ref{p:torusactionlift}.  We therefore concentrate on the $W$-factor.  Since $\rho(\gamma)$ acts algebraically on $W$, if we denote by $g_\gamma : \C^r \to K^\mathbb{C}$ the holomorphic map
\begin{align*}
v \mapsto f(\gamma \acts v) \rho(\gamma) f(v)^{-1} \rho(\gamma)^{-1}
\end{align*}
appearing in the last expression above, it suffices to show that $g_\gamma$ descends to an algebraic map $\tilde{g}_\gamma : V \to K^\C$.  Given that $f : \C^r \to S$ is a homomorphism of complex Lie groups, as noted above, it is straightforward to check that $g_\gamma$ is one as well.  We also observed above that if $v = \lambda \in \Lambda$, then $f(v)^{-1} = \rho(\lambda)$ and so 
\begin{align*}
g_\gamma(\lambda) = f( \gamma \acts \lambda) \rho(\gamma) \rho(\lambda) \rho(\gamma^{-1}) = \rho(\gamma \acts \lambda)^{-1} \rho( \gamma \lambda \gamma^{-1}) = e,
\end{align*}
as $\gamma \lambda \gamma^{-1} = \gamma \acts \lambda$ in $\pi = \Gamma \ltimes \Lambda$.  This shows that $g_\gamma$ descends to a homomorphism of complex Lie groups $\tilde{g}_\gamma : V = \C^r/\Lambda \to K^\C$.  But any such map is automatically algebraic (for example, in the proof of \cite[III.(8.6) Proposition]{BtD}, note that $V = T^\C$ and we may embed $K^\C \hookrightarrow GL_m(\C)$ for some $m$ as an affine subvariety), which proves that the action is algebraic.
\end{proof}

\begin{rmk}
In the introduction to \cite{Totaro2003}, the question is posed of whether a closed manifold admits a good complexification if and only if it admits a Riemannian metric of non-negative curvature.  The above shows that if it admits a \emph{K\"ahler} metric of non-negative holomorphic bisectional curvature, hence non-negative bisectional curvature (Remark \ref{r:secbis}), then a good complexification exists.
\end{rmk}

\subsection{Relation to adapted complex structures}

We have seen that given a K\"ahler manifold, one may consider its canonical complex extension, which is a complex manifold structure on the cotangent bundle.  The K\"ahler metric has an underlying Riemannian metric which, by Proposition \ref{t:anrep}, we may take to be real analytic.  Thus, there is also well-defined adapted complex structure on a neighbourhood of the zero section in the tangent bundle.  Using the metric, we may identify the tangent and cotangent bundles, so when this Grauert tube is in fact entire, we may ask about the relation between these two complex structures on the cotangent bundle.  It is our aim here to make a statement in this direction, beginning with the two following observations.

It is a particular case of \cite[Example 2.1]{LempertSzoke1991} that if $\Lambda$ is a lattice in $\R^n$ with the flat metric, then the adapted complex structure on $T(\R^n/\Lambda)$ is isomorphic to $\C^n/\Lambda$, which by Lemma \ref{l:cstar} is in turn biholomorphic to $(\C^\times)^n$.  Therefore, by Proposition \ref{p:Ztorus}, if $X$ is a flat complex torus of dimension $n$, then the adapted complex structure from the (standard) K\"ahler metric is biholomorphic to the canonical extension.

Consider now the case of a flag manifold $M$. By Section \ref{s:auts}, whose notation we adopt, we may write $M = U/K = G/P$.  \cite[Theorem 2.2]{Szoke1998} states that the adapted complex structure on $TM = T(U/K)$  arising from a $U$-bi-invariant metric on $U$ is biholomorphic to $G/L$. In particular, this holds if $M$ is an irreducible Hermitian symmetric space with the symmetric metric. In this case, by Proposition \ref{p:fvcxfn} and Remark \ref{r:irredhss} the adapted complex structure coincides with the canonical extension as well.

By taking $M$ as in Proposition \ref{p:bisfinitecover}, we would therefore like to conclude that if $(X, \omega)$ is a compact K\"ahler manifold of non-negative bisectional curvature, then the canonical complex extension $Z_{[\omega]}$ is biholomorphic to the adapted complex structure on $TX$. However, in the statement there, which is itself derived from Theorem \ref{t:Mok}, while the metrics on the $\P^N$-factors of $M$ are assumed to be of non-negative bisectional curvature, they are not necessarily invariant.  Without invariance, such a metric in general will not have entire Grauert tube; see \cite[Theorem 2.7, \S4]{Szoke1991} for a discussion of the case $M=S^2 = \P^1$.\footnote{We are grateful to R\'obert Sz\H{o}ke for pointing this out to us.}  However, we are able to prove the following. 

\begin{thm} \label{thm:same_as_adapted} Let $(X, \omega)$ be a compact K\"ahler manifold with non-negative holomorphic bisectional curvature.  Then there is a (real-analytic) K\"ahler form $\widehat\omega$ in the same K\"ahler class as $\omega$, also with non-negative holomorphic bisectional curvature, whose underlying  Riemannian metric admits an entire Grauert tube biholomorphic to the canonical extension $Z_{[\omega]} \cong Z_{[\widehat \omega]}$.

In the decomposition \eqref{e:Mokdecomp}, if one knows a priori either that all the projective space factors carry a (constant) positive multiple of the Fubini--Study metric, or that there are no such factors, then one may take $\widehat \omega = \omega$ above.
\end{thm}
\begin{proof}
For a Riemannian manifold $(Y, g)$, in the following we will often refer to the tangent bundle $TY$ as a complex manifold, with its adapted complex structure implicitly understood (in all cases discussed, $(Y,g)$ will have entire Grauert tube).  If $\Gamma$ is a finite group acting isometrically and freely on $Y$, then, since the adapted complex structure is defined locally, the canonical diffeomorphism
\begin{align} \label{e:acsqt}
(TY)/ \Gamma \xrightarrow{\sim} T(Y/\Gamma)
\end{align}
is in fact a biholomorphism; here, on the right side, we are taking the metric on $Y/\Gamma$ induced from the $\Gamma$-invariant metric $g$. 

Let $(X, \omega)$ be a compact K\"ahler manifold of non-negative bisectional curvature, so that by Theorem \ref{t:Mok}, the universal covering is of the form $\C^k \times M \to X$, where $\C^k$ is endowed with the flat metric.  We will write $M = U/K$ as above and $\pi = \pi_1(X)$.  Assume first that the metric on $M$ arises from a $U$-bi-invariant metric; if one knows a priori that, in the decomposition \eqref{e:Mokdecomp}, $M$ has no projective space factors, or that each factor of $\P^N$ comes with a multiple of the Fubini--Study metric, then this assumption on $M$ is satisfied.  By the observations made in the two paragraphs preceding the theorem, there are biholomorphisms 
\begin{equation} \label{e:idfactors}
Z_{\C^k}  \cong T\C^k \quad  \, \text{ and } \quad\, Z_M  \cong TM,
\end{equation}
from which together with Lemma\ \ref{l:cxfnbas}\eqref{l:cancxfnprod} we conclude that $T(\C^k \times M) \cong Z_{\C^k \times M}$.  The $\pi$-action on $\C^k \times M$ lifts holomorphically to both $T(\C^k \times M)$ (as $\pi$ acts by isometries) and $Z_{\C^k \times M}$ (by Lemma\ \ref{l:cxfnbas}\eqref{l:cxfnpullback}).  Provided that the biholomorphism between these is $\pi$-equivariant, we have
\begin{align*}
TX = T\left( (\C^k \times M)/\pi \right) \cong T(\C^k \times M)/ \pi \cong Z_{\C^k\times M} / \pi \cong Z_{(\C^k \times M)/\pi} = Z_X,
\end{align*}
using \eqref{e:acsqt}, equivariance, and Lemma \ref{l:cxfnbas}\eqref{l:cxfnpullback} for the isomorphisms.  Observe that in this case we have not changed the metric, so this proves the last statement of the Theorem.

To show that the isomorphism $T(\C^k \times M) \cong Z_{\C^k \times M}$ is $\pi$-equivariant, as in the proof of Proposition \ref{p:bisfinitecover}\eqref{p:bisfinitecovera}, we first notice that the $\pi$-action on $\C^k \times M$ is diagonal. The same hence holds for the $\pi$-action on $Z_{\C^k \times M} = Z_\C^k \times Z_M$ stemming from the application of Lemma\ \ref{l:cxfnbas}\eqref{l:cxfnpullback}. Therefore, it suffices to show that each of the two individual identifications in \eqref{e:idfactors} can be taken to be $\pi$-equivariant. Of course, both $T\C^k$ and the canonical complex extension $Z_{\C^k}$ are both biholomorphic to $\C^{2k}$ (recall Section \ref{s:Euctori}).  Again in the proof of Proposition \ref{p:bisfinitecover}, it was noted that the $\pi$-action on $\C^k$ is free, so we can be fairly explicit about the $\pi$-action on $\C^k$, so it is a matter of writing out the maps explicitly to see that the biholomorphism $T \C^k \to Z_{\C^k}$ can be taken to be $\pi$-equivariant; this is done in Lemma \ref{l:toruslift} below.  The canonical complexification of the $M$-factor was computed in Corollary \ref{cor:canonical_complexification_IHSS}, and the (uniquely determined) lift of the action of any group of automorphisms was identified in Sections \ref{subsubsect:flag_automo} and \ref{s:auts} above.  Equivariance of the biholomorphism given in \cite[Theorem 2.2]{Szoke1998} is verified in Lemma \ref{l:Mlift} below.  

To deal with the cases where the projective space factors of $M$ have non-standard metrics, we would simply like to average these metrics to obtain an invariant metric and then proceed as above; to see that this can be done in a matter compatible with the $\pi$-action, we will need to look a little more closely at the form automorphisms of $M$ take.  To this end, we will write
\begin{align*}
(M, \nu) = \prod_{i=1}^p (\P^{N_i}, \theta_i)^{r_i} \times (M', \nu')
\end{align*}
where we write $M' = U'/K'$ as a quotient of compact Lie groups, with the metric $\nu'$ is induced by a $U'$-bi-invariant metric, and $(\P^{N_i}, \theta_i)$ is non-isometric to $(\P^{N_j}, \theta_j)$ if $i \neq j$; note that this does not preclude the possibility that $N_i = N_j$ with $i \neq j$, but if this arises, one has $g^* \theta_i \neq \theta_j$ for all $g \in PGL(N_i+1) = \Aut(\P^{N_i})$.  By the arguments of \cite[\S3]{BeauvilleRmks}, we may decompose
\begin{equation} \label{eq:automogroupdecomposition}
\holIso(M, \nu) = \prod_{i=1}^p \left( \holIso(\P^{N_i}, \theta_i)^{r_i} \rtimes \mathfrak{S}_{r_i} \right) \times \holIso(M', \nu');
\end{equation}
that is, $\holIso(M, \nu)$ acts by biholomorphic isometry in each factor, possibly permuting components in isometric factors.  Of course, one has $\holIso(\P^{N_i}, \theta_i) \leq PGL_{N_i+1}(\C)$.  Let $\pi_i \leq PGL_{N_i+1}(\C)$ be the subgroup of all $\gamma \in PGL_{N_i+1}(\C)$ which appear in some $\holIso(\P^{N_i}, \theta_i)$-component of $\holIso(\P^{N_i}, \theta_i)^{r_i} \rtimes \mathfrak{S}_{r_i}$, when the images of all elements of $\pi$ are decomposed according to \eqref{eq:automogroupdecomposition}.  Since $\pi_i \leq \holIso(\P^{N_i}, \theta_i)$ is contained in a compact group, by conjugating if necessary, we may assume that $\pi_i$ is contained in the maximal compact subgroup $U_i:=PU(N_i+1) = PSU(N_i + 1) \subset PGL_{N_i+1}(\mathbb{C})$, which acts transitively on $\mathbb{P}^{N_i}$. The metric $\tilde{\theta}_i$ we obtain on $\P^{N_i}$ by averaging the metric $\theta_i$ over the compact group $U_i$ must be a scalar multiple of the Fubini--Study metric, and hence arises from a choice of $U_i$-bi-invariant metric on $U_i$.

Now, if we average the (product) metric $\nu$ on $M$ over the group $U = \prod_{i=1}^p U_i^{r_i}$, we obtain a new real-analytic K\"ahler metric $\tilde \nu$ so that isometrically and biholomorphically
\begin{align*}
(M, \tilde{\nu}) = \prod_{i=1}^p (\P^{N_i}, \tilde{\theta}_i)^{r_i} \times (M',\nu').
\end{align*}
The metric $\tilde{\nu}$ is $\pi$-invariant, since it is obviously invariant under any permutation of the factors in each $(\P^{N_i}, \tilde{\theta}_i)^{r_i}$, and since it is invariant in any factor $(\P^{N_i}, \tilde{\theta}_i)$ by construction. Furthermore, as we averaged over the connected compact group $U$, the cohomology class of $\tilde{\nu}$ in $H^1(M, \Omega_M^1) = H^{1,1}(M) \subset H^2(M, \mathbb{C})$ is the same as that of $\nu$ \cite[Lemma 2.3]{Zhang}, so that the induced canonical complexifications are isomorphic as affine bundles over $M$.  Moreover, $\tilde \nu$ obviously has non-negative bisectional curvature. The group $U$ acts transitively on $M$ and the metric $\tilde{\nu}$ arises from a bi-invariant metric on $U$. 

The K\"ahler metric $\widehat \omega$ on $X = (\C^k \times M)/\pi$ induced by the $\pi$-invariant metric $\tilde \nu$ is real-analytic, has non-negative bisectional curvature, so to put ourselves into the situation dealt with in the first paragraphs of the proof, it suffices to show that it belongs to the same cohomology class in $H^{1,1}(X)$ as $\omega$.  For this, we saw above that the $U$-invariant $\tilde{\nu}$ is in the same cohomology class in $H^{1,1}(M)$ as the $\pi$-invariant $\nu$.  If we write $\nu - \tilde{\nu} = \delbar \alpha$ for some $(1,0)$-form $\alpha$, then it suffices to show that $\alpha$ can be taken to be $\pi$-invariant; for then, $\alpha$ would pull back to a $\pi$-invariant primitive on $\C^k \times M$, which would then descend to one on $X$.  Now, if we take $J$ to be the closure of the image of $\pi$ in $U$, then $J$ is a (possibly disconnected) compact Lie subgroup of $U$ and $\nu$ is $J$-invariant (as a simple argument by taking limits shows).  Averaging $\alpha$ over $J$ then yields a $J$- and hence $\pi$-invariant $\tilde{\alpha}$ and one has $\delbar \tilde{\alpha} = \delbar \alpha = \nu - \tilde{\nu}$.
\end{proof}

To complete the proof of Theorem \ref{thm:same_as_adapted} it remains to establish the necessary equivariance properties of the two identifications in \eqref{e:idfactors}.

\begin{lem} \label{l:toruslift}
For $\C^k$ with the flat metric, consider its real tangent bundle $TS$ taken with the adapted complex structure, and its canonical complex extension $Z_{\C^k} \cong \C^{2k}$ as described in Section \ref{s:Euctori}.  For a holomorphic isometry $f : \C^k \to \C^k$, denote by $f^T$ the natural lift of $f$ to $T\C^k$ and let $f^Z$ be the lift of $f$ to $Z_{\C^k}$ as described in the the proof of Proposition \ref{p:torusactionlift}.  There is a biholomorphism $\Phi : TS \xrightarrow{\cong} Z_S$ such that, for any such $f$, the following diagram commutes
\begin{align}\begin{gathered} \label{e:ZTliftdiag}
\xymatrix{ T\C^k \ar[r]^\Phi \ar[d]_{f^T} & Z_{\C^k} \ar[d]^{f^Z} \\ T\C^k \ar[r]_\Phi & Z_{\C^k}. }
\end{gathered}
\end{align}
\end{lem}

\begin{proof}
We write $T\C^k = \C^k \otimes_\R \C$, so that we are realising $\C^k$ as a real vector space, and the adapted complex structure comes from tensoring with $\C$ (cf.\ \cite[Example 2.1]{LempertSzoke1991}).  One thinks of elements of the form $v \otimes 1$ as ``base directions'' and those of the form $v \otimes i$ as ``tangent directions''.  We saw in Section \ref{s:Euctori} that $Z_{\C^k} = \C^{2k}$.  We define $\Phi : TS \to Z_{\C^k} = \C^{2k}$ by
\begin{align*}
v \otimes \alpha \mapsto \left(\alpha v, -\tfrac{i}{2} \alpha \bar{v} \right).
\end{align*}

By the Bieberbach theorem, since one has $U(k) = GL_k(\C) \cap O(2k)$, any holomorphic isometry $f$ of $\C^k$ can be factored into a composition of a translation $\tau$ and a unitary linear map $\sigma$.  We consider these two types of isometries separately.  Suppose first that $\tau$ is translation by $\mu \in \C^k$, i.e.,
\begin{align*}
\tau(z) = z + \mu.
\end{align*}
Then the natural lift $\tau^T$ to $T\C^k$ is
\begin{align*}
v \otimes \alpha \mapsto  v \otimes \alpha + \mu \otimes 1,
\end{align*}
since translation does not affect the tangent directions
The lift $\tau^Z : Z_S \to Z_S$, using the coordinates \eqref{e:Cnholcoords}, is given by
\begin{align*}
(z,w) \mapsto \left( z + \mu, w - \tfrac{i}{2} \bar{\mu} \right),
\end{align*}

If $\sigma : \C^k \to \C^k$ is a unitary linear map, then the lift $\sigma^T$ to $T\C^k$ is given by $v \otimes \alpha \mapsto \sigma(v) \otimes \alpha$, while $\sigma^Z : Z_{\C^k} \to Z_{\C^k}$ is given by $(z,w) \mapsto \left( \sigma(z), \overline{\sigma}(w) \right)$, noting that the second coordinate corresponds to the cotangent fibre directions, and hence the unitary $\sigma$ acts by $(\sigma^t)^{-1} = \overline{\sigma}$.  It is straightforward to check that diagram \eqref{e:ZTliftdiag} commutes in each case.
\end{proof}

In order to check equivariance for the flag manifold factor, let $M = G/P = U/K$ be as in the proof of Theorem \ref{thm:same_as_adapted}, where we assumed that the K\"ahler metric on $M$ comes from a $U$-bi-invariant metric on $U$.  Let us write $M = \prod (M_j, \nu_j)$ as a product of irreducible Hermitian symmetric spaces of compact type, say with $M_j = U_j/K_j$, for some compact Lie groups $K_j \leq U_j$, and the metric $\nu_j$ coming from a $U_j$-bi-invariant metric (we can take $U = \prod U_j$ and $K = \prod K_j$).  Then as observed in Remark \ref{r:irredhss}, $\nu_j$ is a positive scalar multiple of the pullback of the Fubini--Study metric $\omega_j$ under some projective embedding $M_j \hookrightarrow \P^N$, and so by Remark\ \ref{rmk:scalar_multiples_same_bundle} we have \[Z_{[\nu_j]} = Z_{[\omega_j]} \cong U_j^\C / K_j^\C,\] the isomorphism mapping the section coming from the K\"ahler metric to the image of $U/K$ in $U_j^\C / K_j^\C$. Taking products with Lemma \ref{l:cxfnbas}\eqref{l:cancxfnprod}, we obtain $Z_M \cong G/L$ with the section mapping to $U/K \hookrightarrow G/L$, see Proposition \ref{p:fvcxfn}.

By Corollary \ref{c:Mactionlift} and Lemma \ref{l:uniquelift}, $\Aut M = E_P \ltimes G$ uniquely lifts to an action on $Z_M$; hence so does the action of $\holIso(M)$, which is contained in $E_P \ltimes U$ by Proposition \ref{p:holisoM}.  On the other hand, $\Aut M$ also has a natural lift to an action on $TM$; since $\holIso(M) \leq \Aut M$ acts by isometries, this action on $TM$ is by biholomorphism with respect to the adapted complex structure.

\begin{lem} \label{l:Mlift}
In this setup, the biholomorphism $\Phi : TM \to Z_M$ of \cite[Theorem 2.1]{Szoke1998} intertwines the two lifted actions discussed above.  In particular, if $\Gamma$ is any group acting on $M$ via biholomorphic isometries, then $\Phi$ intertwines the lifted $\Gamma$-actions on $TM$ and $Z_M$. 
\end{lem}

\begin{proof}
Following \cite[\S2]{Szoke1998}, though changing the notation, let us write $\mathfrak{u}$ and $\mathfrak{k}$ for the Lie algebras of $U$ and $K$, respectively, and write 
\begin{align*}
\mathfrak{u} = \mathfrak{k} \oplus \mathfrak{m}
\end{align*}
for some $K$-invariant subspace $\mathfrak{m} \subseteq \mathfrak{u}$.  Then we will realise $TM$ as the associated vector bundle $TM = U \times^K \mathfrak{m}$,
and the biholomorphism $\Phi : TM \to U^\C/K^\C = G/L$ of \cite[Theorem 2.1]{Szoke1998} is given by
\begin{align*}
[a, \xi] \mapsto a \exp(i \xi) L.
\end{align*}
As mentioned, we identify $Z_M$ with the target of $\Phi$ and it is via this identification that we will realise the lifted action of Corollary \ref{c:Mactionlift}.

From Proposition \ref{p:holisoM}, one has $\holIso(M) \leq E_P \ltimes U$; it therefore suffices to show that $\Phi$ intertwines the actions of both elements of $U$ as well as those of $E_P$.  An element $a \in U$ acts on $M = U/K$ and $Z_M = G/L$ via left multiplication $a \acts uL = auL$. The lifted action on $TM$ is $a \acts [u, \xi] = [au, \xi]$. For this type of automorphism, it is obvious that $\Phi$ intertwines the action. The second type arises from $\phi \in E_P$, which has a realisation as an automorphism $\phi : U \to U$, see discussion preceding Proposition \ref{p:holisoM}.  Then 
\begin{align*}
\Phi\left( \phi \acts [u, \xi] \right) = \phi(u) \exp \left( i d\phi(\xi) \right) L = \phi(u) \phi \left( \exp(i \xi) \right) L = \phi\left( u \exp(i\xi) \right)L = \phi \acts \Phi([u,\xi]). & \qedhere
\end{align*}
\end{proof}

\section{Big tangent bundles and affine canonical extensions}

In this section we start investigating the question of necessary assumptions on a manifold to have a canonical extension with many holomorphic or regular functions. Moreover, as a by-product of our investigation, we give a new proof of a result of Yang \cite[Theorem 4.5]{Yang2017}.

\subsection{Affineness of the canonical extension and bigness of the tangent bundle} \label{subsect:big_tangent}

We now consider the case where $X$ is a smooth complex projective variety.  In this case, it is possible to show that if $Z_X$ is affine, then the tangent bundle $\Theta_X$ must be big.  In particular, this applies to flag varieties as discussed in Section \ref{s:fvHss}.

In order to obtain the result just mentioned, we give a more general statement about affine open subvarieties in projective varieties, which in turn makes essential use of the following result of Goodman \cite[\S I Theorem 1]{Goodman1969}, which is also stated and proved at \cite[\S II.6 Theorem 6.1]{Hartshorne}.

\begin{thm} \label{t:Goodman} 
Let $Y$ be a closed subset of a complete integral scheme $X$ and let $U = X - Y$.  Then $U$ is affine if and only if there exists a closed subscheme $Z \subseteq Y$ such that $\overline{Y} = f^{-1}(Y)$ is the support of an effective ample divisor on $\overline{X}$, where $f : \overline{X} \to X$ is the blow-up of $X$ with centre $Z$.
\end{thm}

\begin{prop}
Let $Y$ be a normal projective variety and $D$ a smooth connected divisor on $Y$.  If $Y \setminus D$ is an affine variety, then the normal bundle $\mathcal{N}_{D/Y} \cong \O_D(D) = \O_Y(D)|_D$ is big.
\end{prop}

\begin{proof}
By Theorem \ref{t:Goodman}, there exists a birational modification $f : \overline{Y} \to Y$ of $Y$ with centre $Z$ contained in $D$ and an ample effective divisor $A$ on $\overline{Y}$ with support $f^{-1}(D)$. By blowing-up further if necessary, we may assume that $\overline{Y}$ is smooth; moreover, we may assume that $\codim_Y Z \geq 2$.

By taking a sufficiently large multiple, we may assume that $A$ is very ample and that it takes the form $A = k f^*D -N +P$, 
where $k \in \Z_{>0}$ and $N$ and $P$ are effective divisors which are exceptional for $f$.  Thus,
\begin{align*}
f_* \O_{\overline{Y}}(A) = f_* \O_{\overline{Y}}( kf^* D - N + P) = \O_Y(kD) \otimes_{\O_Y} f_* \O_{\overline{Y}}( -N + P).
\end{align*}
But now, $\O_{\overline{Y}}(-N) \subseteq \O_{\overline{Y}}$, so by left exactness of $f_*$, one has $f_* \O_{\overline{Y}} (-N+P) \subseteq f_* \O_{\overline{Y}}(P)$. We then claim that $f_* \O_{\overline{Y}}(P) = \O_Y$.  This statement is a consequence of a special case of \cite[Lemma 1-3-2]{KMM1987}, however, one can also give the following simple argument.  First, $\O_{\overline{Y}}(P)$ may be realised as the subsheaf of the sheaf of rational functions on $\overline{Y}$ with poles bounded by $P$.  However, as $f$ is birational and $f(P) \subseteq Z$, when we realise these as rational functions on $Y$ via $f$, they are well-defined except possibly over $Z$, which is of codimension $\geq 2$ in the normal variety $Y$, so  extend to functions on $Y$ without any poles, thus yielding sections of $\O_Y$.

It follows that $f_* \O_{\overline{Y}}(A) \subseteq \O_Y(kD)$ and hence we have an inclusion of global sections
\begin{align} \label{e:globsecinc}
V := H^0(\overline{Y}, A) \hookrightarrow H^0(Y, \O_Y(kD)).
\end{align}
Since $A$ is very ample, its global sections yield a projective embedding $\phi_A : \overline{Y} \hookrightarrow \P^N$; but considering the linear system $|V|$ on $Y$ induced by the inclusion \eqref{e:globsecinc}, we obtain a rational map $\phi_{|V|} : Y \dashrightarrow \P^N$ yielding a factorisation $\phi_{kD}$ of $\phi_A$ that fits into a commutative diagram
\begin{align*}
\xymatrix{
\overline{D} \ar[d]_f \ar[r]^{\overline{\iota}} & \overline{Y} \ar[d]_f \ar[r]^{\phi_A} & \P^N \\
D \ar[r]^\iota & Y \ar@{-->}[ur]_{\phi_{|V|}} &
}
\end{align*}
where $\overline{D} \subseteq \overline{Y}$ is the strict transform of $D$ and $\iota$ and $\overline{\iota}$ are the respective inclusion maps.  Since $f$ is an isomorphism away from $f^{-1}(Z)$, the diagram already shows that the indeterminacy locus of $\phi_{V}$ is contained in $Z$.  In particular, since $\phi_A \circ \overline{\iota}$ is an embedding, $\iota \circ \phi_{|V|}$ is an embedding at least on $D \setminus Z$.  But $\iota \circ \phi_{|V|}$ is the rational map obtained from the restriction of the linear system $|V|$ to $D$, which is a linear subsystem of $\O_Y(kD)|_D = \mathcal{N}_{D/Y}^k$; this says that  $\mathcal{N}_{D/Y}$ is big.
\end{proof}

\begin{cor} \label{c:big} Suppose $X$ is a smooth projective variety such that $Z_X$ is affine.  Then $\Theta_X$ is big.  In particular, if $X$ is a flag variety, then $\Theta_X$ is big.
\end{cor}

\begin{proof}
In the situation of Section \ref{ss:cancxfn}, we take $Y = \P W$ and $D = \P \Omega_X^1$.  Then by Lemma \ref{l:extlem}(\ref{l:PWminusPV}), $Z_X \cong Y \setminus D$ and $\O_{\P W}(D) = \O_{\P W}(1)$.  Of course, $\O_{\P W}(1)|_{\P \Omega_X^1} = \O_{\P \Omega_X^1}(1) = \mathcal{N}_{\P \Omega_X^1/ \P W}$.  The statement that $\O_{\P \Omega_X^1}(1)$ is big is the definition of $\Theta_X$ being big. The last statement comes from Proposition \ref{p:fvcxfn}.
\end{proof}

\begin{rmk}
Our last statement is the second part of \cite[Corollary 1.3]{Hsiao2015}, which is obtained via much different methods.  The first part of Hsiao's statement suggests the question of whether $Z_X$ is affine when $X$ is a toric variety.
\end{rmk}

The following gives a first indication of the place of varieties with big tangent bundles in the classification theory of higher-dimensional projective varieties.

\begin{prop} \label{p:Tbiguniruled}
If $X$ is a projective manifold with the property that $\Theta_X$ is big, then $X$ is uniruled.  In particular, $X$ having big tangent bundle implies that its Kodaira dimension is negative. 
\end{prop}

\begin{proof}
Assume that $X$ is not uniruled and fix an ample line bundle $A$ on $X$. By Miyaoka's generic semipositivity theorem (see, e.g., \cite[Remark 6.3.34]{Lazarsfeld}), for a sufficiently general curve $C \subseteq X$ arising as a complete intersection of $n-1$ divisors in $|rA|$ for $r \gg 1$, the restriction $\Omega_X^1|_C$ is nef.  But then $(\Omega_X^1|_C)^{\otimes m} = (\Omega_X^1)^{\otimes m}|_C$ is nef for any $m > 0$ by \cite[Theorem 6.2.12(iv)]{Lazarsfeld}, and hence any quotient of $(\Omega_X^1)^{\otimes m}|_C$ is also nef by \cite[Theorem 6.2.12(i)]{Lazarsfeld}.

On the other hand, since $\Theta_X$ is big, \cite[Example 6.1.23]{Lazarsfeld} says that for sufficiently large $m$, we have $H^0( X, \Sym^m \Theta_X \otimes A^\vee) \neq 0$. Thus, $A$ is a subsheaf of $\Sym^m \Theta_X$.  Furthermore, we may identify the sub-bundle of symmetric $m$-tensors in $\Theta_X^{\otimes m}$ with $\Sym^m \Theta_X$ (this works in characteristic $0$).  Therefore, $A$ is a subsheaf of $\Theta_X^{\otimes m}$; dualising, we hence obtain an anti-ample quotient of $(\Omega_X^1)^{\otimes m}$ and restricting to $C$ above, we obtain a contradiction.
\end{proof}

\begin{rmk}
In the statement of Proposition \ref{p:Tbiguniruled}, we have assumed that $X$ is projective; however, with the assumption that $\Theta_X$ is big, this is automatic.  Assuming only that $X$ is compact K\"ahler, the same is true for the projective bundle $\P \Omega_X^1$.  The existence of a big line bundle, namely $\O_{\P \Omega_X^1}(1)$, implies that $\P \Omega_X^1$ is a Moishezon manifold.  Now, the image of a Moishezon manifold under a holomorphic map is Moishezon; hence, $X$ is compact K\"ahler and Moishezon, and therefore projective by Moishezon's theorem. 
\end{rmk}

\subsection{Non-bigness of tangent bundles of fibre bundles over tori} 

Here, we present a somewhat more elementary proof of the result of Yang \cite[Theorem 4.5]{Yang2017} stating that a compact K\"ahler manifold of non-negative bisectional curvature with big tangent bundle is a product of Hermitian symmetric spaces of compact type; see Corollary\ \ref{cor:yang} below for the precise formulation. As in Theorems \ref{t:bisholStein} and \ref{thm:good_complexification}, we will again apply the relevant structure theory as provided by Proposition \ref{p:bisfinitecover}.

\begin{prop} \label{p:XTMnotbig}
 Let $T = \C^k / \Lambda$ be a complex torus of dimension $k > 0$, so that we have a natural identification $\pi_1(T) = \Lambda$.  Let $M$ be any compact complex manifold on which $\Lambda$ acts by biholomorphism via a compact subgroup; i.e., the homomorphism $\Lambda \to \Aut M$ factors through a compact subgroup of $\Aut M$.  Let
\begin{align*}
X := \C^k \times^\Lambda M = (\C^k \times M)/ \Lambda.
\end{align*}
Then $\Theta_X$ is not big.
\end{prop}

\begin{cor}[\protect{\cite[Thm.~4.5]{Yang2017}}] \label{cor:yang}
If $X$ is a compact K\"ahler manifold of non-negative holomorphic bisectional curvature with big tangent bundle, then $X$ is a product 
\begin{align*}
(\P^{N_1}, \omega_1) \times \cdots \times (\P^{N_\ell}, \omega_\ell) \times (M_1, \eta_1) \times \cdots \times (M_k, \eta_k)
\end{align*}
where, for $1 \leq i \leq \ell$, $\omega_i$ is a K\"ahler metric of non-negative bisectional curvature on $\P^{N_i}$,  and where each $(M_j, \eta_j)$ is an irreducible Hermitian symmetric space of compact type.
\end{cor}

\begin{proof}
If $X$ is of non-negative bisectional curvature, then by Proposition \ref{p:bisfinitecover}\eqref{p:bisfinitecoverT}, it has a finite Galois cover $Y$ which is of precisely the form as described in Proposition \ref{p:XTMnotbig} above.  If $k > 0$, then the Proposition and Lemma \ref{l:etbig} yield a contradiction to the assumption that $\Theta_X$ is big.  Therefore, we must have $k = 0$ and so $Y$ is a product as in the statement, and $X$ is a finite quotient of $Y$; but by Proposition \ref{p:Hssbas}, this must be a quotient by the trivial group.
\end{proof}

The rest of this section will be devoted to providing a proof of Proposition \ref{p:XTMnotbig}.  We will write $\dim M = n-k$ so that $\dim X = n$.  Let $f : X \to T$ denote the projection.  Then since $\Theta_T \cong \O_T^{\oplus k}$, one has $f^* \Theta_T \cong \O_X^{\oplus k}$, and the relative tangent sequence associated to $f$ is of the form
\begin{align} \label{e:XTMreltan}
0 \longrightarrow \Theta_{X/T} \longrightarrow \Theta_X \longrightarrow \O_X^{\oplus k} \longrightarrow 0.
\end{align}

Let $P_X$ and $P_{X/T} \in \mathbb{Q}[z]$ be the relevant Hilbert functions, i.e.,
\begin{align*}
P_X(m) & := h^0(X, \Sym^m \Theta_X) = h^0( \P \Omega_X^1, \O_{\P \Omega_X^1}(m) ) \\
P_{X/T}(m)& := h^0(X, \Sym^m \Theta_{X/T}) = h^0( \P \Omega_{X/T}^1, \O_{\P \Omega_{X/T}^1}(m) ). 
\end{align*}

Recall that $P_X(m) = O(m^{2n-1})$ and $P_{X/T}(m) = O(m^{2n - k - 1})$ as $\P \Omega_X^1$ and $\P \Omega_{X/T}^1$ are of dimensions $2n-1$ and $2n-k-1$, respectively.  Furthermore, each of the bundles being big is equivalent to the respective limits superior,
\begin{align*}
& \limsup_{m \to \infty} \frac{P_X(m)}{m^{2n-1}} & &\text{and}&& \limsup_{m \to \infty} \frac{ P_{X/T}(m) }{ m^{2n-k-1} }
\end{align*}
being strictly positive (cf.\ Definition \ref{d:big}).

Proposition \ref{p:XTMnotbig} is proved by verifying the following two lemmata.

\begin{lem}
If $\Theta_X$ were big, then $\Theta_{X/T}$ would also be big.
\end{lem}

\begin{proof}We begin with the following statement from multilinear algebra \cite[Chap.\ II, Ex.\ 5.16]{Hartshorne1977}. Note that our indexing differs from that given in the reference.

\bigskip

\noindent
\emph{Claim.} 
Let $Y$ be any complex manifold and let 
\begin{align*}
0 \longrightarrow U \longrightarrow V \longrightarrow W \longrightarrow 0
\end{align*}
be an exact sequence of (holomorphic) vector bundles on $Y$.  Then for any $m \in \Z_{> 0}$, there is a filtration
\begin{align*}
0 = F_{-1} \subseteq F_0 \subseteq F_1 \subseteq \cdots \subseteq F_m = \Sym^m V
\end{align*}
by sub-bundles such that $F_i/F_{i-1} \cong \Sym^{m-i} U \otimes \Sym^i W$.

\bigskip

\noindent
The claim is proved by taking $F_i$ to be the image of 
\begin{align*}
\bigoplus_{j=0}^i U^{\otimes m-j} \otimes V^{\otimes j} = U^{\otimes m} \oplus U^{\otimes m-1} \otimes V \oplus \cdots \oplus U^{\otimes m-i} \otimes V^{\otimes i} \subseteq V^{\otimes m}
\end{align*}
in $\Sym^m V$ under the symmetrisation map.

\bigskip

Applying the Claim to the exact sequence \eqref{e:XTMreltan}, we obtain
\begin{align} \label{e:grTheta}
F_i/F_{i-1} \cong \Sym^{m-i} \Theta_{X/T} \otimes \Sym^i \O_X^{\oplus k} \cong \left( \Sym^{m-i} \Theta_{X/T} \right)^{\oplus \binom{k+i-1}{i}}.
\end{align}

We have 
\begin{align*}
P_X(m) & = h^0(X, \Sym^m \Theta_X) = h^0(X, F_m) \leq h^0(X, F_m/F_{m-1}) + h^0(X, F_{m-1}) \\
& = \binom{m+k-1}{m} + h^0(X, F_{m-1}),
\end{align*}
and a straightforward induction, repeatedly using \eqref{e:grTheta}, shows that 
\begin{align} \label{e:PXineq}
P_X(m) & \leq \binom{m+k-1}{k-1} + \binom{m+k-2}{k-1} P_{X/T}(1) + \binom{m+k-3}{k-1} P_{X/T}(2) + \cdots \nonumber \\
& \quad \quad + \binom{k+1}{k-1} P_{X/T}(m-2) + k P_{X/T}(m-1) + P_{X/T}(m) \nonumber \\
& = \binom{m+k-1}{k-1} + \sum_{i=1}^m \binom{m + k - i - 1 }{ k - 1 } P_{X/T}(i).
\end{align} 

Let us assume that $\Theta_{X/T}$ is not big, i.e., $\lim \sup_m P_{X/T}(m)/ m^{2n-k-1} = 0$, so that we want to prove that $\lim \sup_m P_X(m)/m^{2n-1} = 0$.  It therefore suffices to divide \eqref{e:PXineq} by $m^{2n-1}$ and then show that the limit is zero.  Observe that for each $i \geq 0$, $\binom{m+k-i-1}{k-1}$ is a polynomial in $m$ of degree $k-1$, and in fact, the leading term is $m^{k-1}/(k-1)!$ in each case.  Thus, the first term, when divided by $m^{2n-1}$, will vanish in the limit, so it suffices to look at the sum.  Dividing by $m^{2n-1}$, we obtain
\begin{align*}
\sum_{i=1}^m \frac{1}{m^k} \binom{m + k - i - 1 }{ k - 1 } \frac{P_{X/T}(i)}{ m^{2n-k-1} }.
\end{align*}
Since the $P_{X/T}(i)/ m^{2n-k-1}$ can be made uniformly small by taking $m$ large enough (i.e., independently of $i$), it suffices to show that $\lim_m \sum_{i=1}^m \binom{ m + k - i - 1}{ k-1 }/ m^k$ is bounded, but this again follows from the fact that the given binomial coefficients are polynomials of degree $k-1$.
\end{proof}

\begin{lem} \label{lem:XTnotbig}
$\Theta_{X/T}$ is not big.
\end{lem}

\begin{proof}
The $\Lambda$-action on $M$ induces one on the tangent bundle $\Theta_M$ (as a complex manifold) and hence on $\Sym^m \Theta_M$ for any $m \in \Z_{>0}$; of course, the action is linear on fibres, and the projections are equivariant for the action on $M$.  We have $\Theta_{(\C^k \times M)/\C^k} = p_M^* \Theta_M$, and hence since $\Lambda$ is a discrete group
\begin{align*}
\Theta_{X/T} = (\C^k \times \Theta_M)/\Lambda = \C^k \times^\Lambda \Theta_M
\end{align*}
and of course the projection to $X$ is obtained by applying the functor $\C^k \times^\Lambda -$ to $\Theta_M \to M$.  We can repeat this for all symmetric products:  for all $m \in \Z_{>0}$, one has
\begin{align*}
\Sym^m \Theta_{X/T} \cong \C^k \times^\Lambda \Sym^m \Theta_M.
\end{align*}

Now, let $E$ be any holomorphic vector bundle over $M$ on which there is a $\Lambda$-action, also acting via a compact group, which is linear on the fibres and for which the projection to $M$ is equivariant.  Of course, over $\Lambda$-invariant open sets of $M$, there is a $\Lambda$-action on the sections of $E$ by $(\lambda \acts s)(m) = s \left( (-\lambda) \acts m \right)$,
where $\lambda \in \Lambda$, $m \in M$, and $s$ is a section of $E$; we are using additive notation for $\Lambda$.

We let
\begin{align*}
E_X := \C^k \times^\Lambda E,
\end{align*}
so that the natural projection to $X$, as before obtained by applying $\C^k \times^\Lambda -$ to the projection $E \to M$, makes it into a vector bundle over $X$. 
In this setup, we have the following general result.
\begin{prop} \label{p:ELam}
There is a natural isomorphism
\begin{align*}
H^0(M, E)^\Lambda = H^0(X, E_X),
\end{align*}
where the left hand side denotes the $\Lambda$-invariant global sections of $E$.
\end{prop}

\begin{proof}[Proof of Proposition~\ref{p:ELam}]
We denote the class of $(u, m) \in \C^k \times M$ modulo $\Lambda$ by $[u, m] \in X$, and similarly for elements of $E_X = \C^k \times^\Lambda E$.  Given $s \in H^0(M, E)^\Lambda$, we define $\sigma_s \in H^0(X, E_X)$ by
\begin{align*}
\sigma_s([u,m]) = [u, s(m)].
\end{align*}
That this is well-defined comes directly from the fact that $s$ is $\Lambda$-invariant, i.e., $s(\lambda \cdot m) = \lambda \cdot s(m)$ for all $\lambda \in \Lambda$, $m \in M$.

Now, consider a global section of $E_X$.  We may pull this back through the projection $\C^k \times M \to X$ and so we get a section of the vector bundle $\C^k \times E \to \C^k \times M$, which must be compatible with the $\Lambda$-action; hence a section of $E_X$ may be identified with a $\Lambda$-equivariant function $\C^k \to H^0(M, E)$.  Since $\Lambda$ is abelian, acting via a compact group, and $H^0(M, E)$ is a finite-dimensional vector space, it decomposes into a sum of $1$-dimensional $\Lambda$-modules; furthermore, since $\Lambda$ acts via a compact subgroup, it must act via $S^1$ on each of these $1$-dimensional modules.  Therefore, a section of $E_X$ may be identified with a tuple of sections of $1$-dimensional local systems on $T$, each arising from a representation of $\Lambda$ in $S^1$.  Each such section corresponds to a holomorphic function on $\C^k$ that is invariant under the $\Lambda$-action. As $\Lambda$ acts via $S^1$ on $\mathbb{C}$, the function is in fact bounded, and hence constant by Liouville's Theorem; but for a constant function to be $\Lambda$-invariant, the representation would have to be trivial. In summary, we can identify the sections of $E_X$ with constant maps $\C^k \to H^0(M, E)^\Lambda \subseteq H^0(M, E)$.
\end{proof}

To complete the proof of Lemma \ref{lem:XTnotbig}, we will take $E = \Sym^m \Theta_M$, so that $E_X = \Sym^m \Theta_{X/T}$.  As above, we write $P_M \in \mathbb{Q}[z]$ for the Hilbert function
\begin{align*}
P_M(m) := h^0(M, \Sym^m \Theta_M) = h^0 \big( \P \Omega_M, \O_{\P \Omega_M}(m) \big).
\end{align*}
As $\dim \P \Omega_M = 2 \dim M - 1 = 2(n-k) -1$, $P_M(m)/ m^{2(n-k)-1}$ is bounded and has a positive limit superior as $m \to \infty$ if and only if $\Theta_M$ is big (which is the case when $M$ is a flag variety, by Corollary \ref{c:big}, the case of interest to us).  Proposition \ref{p:ELam} then states that
\begin{align*}
P_{X/T}(m) = h^0( X, \Sym^m \Theta_{X/T}) \leq h^0(M, \Sym^m \Theta_M) = P_M(m)
\end{align*}
and so
\begin{align*}
\frac{P_{X/T}(m)}{m^{2n-k-1}} \leq \frac{P_M(m)}{m^{2n-k-1}} = \frac{P_M(m)}{m^{2(n-k)-1}} \cdot \frac{1}{m^k}\; \overset{m\to \infty}{\longrightarrow}\; 0
\end{align*}
by our observation above that the first factor in the last expression is bounded.  Hence $\Theta_{X/T}$ cannot be big. 
\end{proof}
This concludes the proof of Proposition~\ref{p:XTMnotbig}.

\vspace{0.6cm}

\appendix

\section{Proof of Proposition \ref{t:anrep}} \label{app:anrep}

We understand that a proof of Proposition \ref{t:anrep} is possible using methods involving the K\"ahler--Ricci flow \cite{B}, however, we present here a more elementary, sheaf-theoretic proof.  A proof similar to ours appears at \cite[Proposition 2.1]{Lempert2017}; we thank Bo Berndtsson for making us aware of this reference.

Let $X$ be a complex manifold of (complex) dimension $n$.  Consider the map of sheaves $\Re : \O_X \to C_X^\infty$ which takes the real part of a holomorphic function.  This sits in an exact sequence
\begin{align*}
\xymatrix{
0 \ar[r] & i \R_X \ar[r] & \O_X \ar[r]^{\Re} & C_X^\infty,
}
\end{align*}
where $\R_X$ is the sheaf of locally constant $\R$-valued functions on $X$.  We set $\K_\infty^1 := \coker \Re$, so that we have a short exact sequence
\begin{align*}
\xymatrix{
0 \ar[r] & \O_X / i\R_X \ar[r] & C_X^\infty \ar[r] & \K_\infty^1 \ar[r] & 0.}
\end{align*}

Replacing $C_X^\infty$ by $C_X^\omega$, we get an analogous short exact sequence and the inclusion $C_X^\omega \hookrightarrow C_X^\infty$ yields a morphism of short exact sequences
\begin{align} \label{E}
\vcenter{ \xymatrix{
0 \ar[r] & \O_X/ i\R_X \ar[r] \ar@{=}[d] & C_X^\omega \ar[r] \ar[d] & \K_\omega^1 \ar[d] \ar[r] & 0 \\
0 \ar[r] & \O_X / i\R_X \ar[r] & C_X^\infty \ar[r] & \K_\infty^1 \ar[r] & 0.
} }
\end{align}
The reason for introducing these objects is as follows.  Suppose that $X$ admits a K\"ahler form $\omega$; then $\omega$ admits local K\"ahler potentials:  there exists an open cover $\{ U_j \}$ of $X$ and functions $\varphi_j \in C^\infty(U_j)$ with $\omega|_{U_j} = i \partial \dbar \varphi_j$
and, of course, the $\varphi_j$ will need to be strictly plurisubharmonic.  On the overlaps, one has $i \partial \dbar (\varphi_j|_{U_j \cap U_k} - \varphi_k|_{U_j \cap U_k}) = 0$, so the $\varphi_j|_{U_j \cap U_k} - \varphi_k|_{U_j \cap U_k}$ are pluriharmonic functions, i.e., precisely those in the image of $\O_X/ i \R_X \to C_X^\infty$.  We may therefore realise the K\"ahler form $\omega$ as a global section of the sheaf $\K_\infty^1$ which is locally represented by strictly plurisubharmonic functions.  When we do this, then the map $C_X^\infty \to \K_\infty^1$ should be thought of as the operator $i \partial \dbar$.
 
To proceed with the proof, we begin by assembling the approximation result we need to find analytic functions ``close enough'' to a given strictly plurisubharmonic function.  This begins with some local arguments.

\begin{lem}
Let $B \subseteq \C^n$ be a ball centred at the origin and let $z_1, \ldots, z_n$ be the standard complex coordinates on $\C^n$.  Suppose $u_1, \ldots, u_{2n}$ is any set of real coordinates on $B$.  Then for any $\epsilon > 0$, there exists $\delta > 0$ such that if $\psi \in C^2(B)$ satisfies
\begin{align} \label{A}
|\psi|, \left| \frac{\partial \psi}{\partial u_\alpha} \right|, \left| \frac{\partial^2 \psi}{\partial u_\alpha \partial u_\beta} \right| < \delta
\end{align}
then
\begin{align} \label{B}
|\psi|, \left| \frac{\partial \psi}{\partial z_\lambda} \right|, \left| \frac{\partial^2 \psi}{\partial z_\lambda \partial z_\mu} \right|, \left| \frac{\partial^2 \psi}{\partial z_\lambda \partial \bar{z}_\mu} \right| < \epsilon
\end{align}
on some possibly smaller ball $B' \subseteq B$ containing $0$.
\end{lem}

For the proof, one may restrict to some closed ball $\overline{B}' \subseteq B$ on which all of
\begin{align*}
\frac{\partial u_\alpha}{\partial z_\lambda}, \frac{\partial u_\alpha}{\partial \bar{z}_\lambda}, \frac{\partial^2 u_\alpha}{\partial z_\lambda \partial z_\mu}, \frac{\partial^2 u_\alpha}{\partial z_\lambda \partial \bar{z}_\mu}
\end{align*}
are bounded.  Since $\partial \psi/ \partial z_\lambda$, etc., can be written in terms of these and the $\partial \psi/ \partial u_\alpha$, etc., one may find such a $\delta$.

Consider now a $C^\omega$ (respectively, $C^\infty$) embedding of the complex manifold $X$ into an open set $\Omega \subseteq \R^k$, for some $k > 2n$, with coordinates $u_1, \ldots, u_k$ on $\R^k$.  Write $\kk := \{ 1, \ldots, k \}$.  Given any $x \in X$, there exists at least one subset $\ii =  \{ i_1, \ldots, i_{2n} \}$ such that $u_\ii := (u_{i_1}, \ldots, u_{i_{2n}})$ is a set of $C^\omega$ (respectively, $C^\infty$) coordinates on $X$ in a neighbourhood of $x$.  By taking a minimum over all such $\ii$, we can state the following.

\begin{lem}
Consider $X \hookrightarrow \Omega$ embedded as above, fix $x \in X$ and choose holomorphic coordinates $z_1, \ldots, z_n$ on $X$ near $x$.  Then given any $\epsilon > 0$, there exists $\delta > 0$ such that if $\psi \in C^2(B)$ on some neighbourhood $B$ of $x$ (in $X$) satisfies \eqref{A} for all $\alpha, \beta \in \ii$, where $\ii \subseteq \kk$ is any subset such that $u_\ii$ is a set of (real) coordinates on $X$ near $x$, then \eqref{B} holds on some neighbourhood $B'$ of $x$ (in $X$).
\end{lem}

Fix $x \in X$ and suppose $\ii \subseteq \kk$ is such that $u_\ii$ gives a set of coordinates on $X$ near $x$.  For $\beta \in \kk \setminus \ii$, then in such a neighbourhood on $X$, there exist uniquely determined $C^\omega$ (respectively $C^\infty$) functions $\tau_\beta$ such that $u_\beta = \tau_\beta (u_\ii)$. If $h \in C^2(\Omega)$, then setting $\psi := h|_X$, for all $\alpha \in \ii$, we have
\begin{align*}
\frac{\partial \psi}{\partial u_\alpha} = \frac{\partial h}{\partial u_\alpha} + \sum_{\beta \in \kk \setminus \ii} \frac{\partial \tau_\beta}{\partial u_\alpha} \frac{\partial h}{\partial u_\beta}
\end{align*}
and we can find similar expressions for the second derivatives.  Thus the expressions for $\partial \psi/ \partial u_\alpha$, etc. are functions of $\partial h/ \partial u_\gamma$, $\gamma \in \kk$.  We can thus obtain the following.

\begin{lem}
Consider $X \hookrightarrow \Omega$ embedded as above, fix $x \in X$ and choose holomorphic coordinates $z_1, \ldots, z_n$ on $X$ near $x$.  Then given any $\epsilon > 0$, there exists $\eta > 0$ such that if $h \in C^2(\Omega)$ satisfies \eqref{A} for all $\alpha, \beta \in \kk$, then \eqref{B} holds on some neighbourhood $B'$ of $x$ (in $X$).
\end{lem}

We will need the following basic local deformation result for strictly plurisubharmonic functions.

\begin{lem}[\protect{\cite[Corollary IX.B.3]{GunningRossi}}] \label{t:GR}
Let $\varphi$ be a strictly plurisubharmonic function defined on a neighbourhood of the origin $0 \in \C^n$.  For some $r > 0$, if $\psi \in C^2( B_r)$ ($B_r$ is an open ball of radius $r$ centred at the origin) is such that $|\psi|$, $|\psi_i|$, $|\psi_{ij}|$, $|\psi_{i \bar{j}}|$ are sufficiently small, then $\varphi - \psi$ is also strictly plurisubharmonic.
\end{lem}

With the observation that the notion of a strictly plurisubharmonic function on a complex manifold is independent of coordinates, by applying Lemma~\ref{t:GR}, we can state the following.

\begin{lem}
Consider $X \hookrightarrow \Omega$ embedded as above, fix $x \in X$ and suppose $\varphi_0$ is some strictly plurisubharmonic function defined in a neighbourhood of $x$ in $X$.  Then there exists $\eta = \eta_x > 0$ such that if $h \in C^2(\Omega)$ satisfies
\begin{align} \label{C}
|h|, \left| \frac{\partial h}{\partial u_\alpha} \right|, \left| \frac{\partial^2 h}{\partial u_\alpha \partial u_\beta} \right| < \eta_x
\end{align}
for all $\alpha, \beta \in \kk$, then $\varphi_0 - h|_X$ is strictly plurisubharmonic on a neighbourhood of $x$ (in $X$).
\end{lem}

Now, suppose we are given $\omega_0 \in H^0(X, \K_\infty^1)$ which is represented everywhere by spsh functions. Then we can find a locally finite, countable open cover $\{ U_m \}_{m \in \N}$ of $X$, with $U_m$ relatively compact in $X$ and $\varphi_{0,m}$ strictly plurisubharmonic on $U_m$ and find $\eta_m > 0$ (without loss of generality, $\eta_m < \tfrac{1}{2}$) such that if $h \in C^2(\Omega)$ satisfies \eqref{C}, then $\varphi_{0,m} - h|_{U_m}$ is also strictly plurisubharmonic on $U_m$.

\begin{lem} \label{t:eta}
There exists $\tilde{\eta} \in C^0(X)$ such that for every $x \in X$, there exists some $m \in \N$ with $x \in U_m$ and $0 < \tilde{\eta}(x) < \eta_m$.
\end{lem} 

The statement we need to prove Theorem \ref{t:anrep} will follow from Whitney's approximation theorem.

\begin{thm}[Whitney's approximation theorem, \protect{\cite[1.6.5]{Narasimhan}}] Let $\Omega \subseteq \R^n$ be open and $f \in C^k(\Omega)$, $0 \leq k \leq \infty$.  Let $\eta \in C^0(\Omega)$ be such that $\eta(x) > 0$ for all $x \in \Omega$.  Then there exists $g \in C^\omega(\Omega)$ such that
\begin{align*}
| D^\alpha f(x) - D^\alpha g(x)| < \eta(x) \quad \textnormal{for } 0 \leq |\alpha| \leq \min \left\{ k, \frac{1}{\eta(x)} \right\}. 
\end{align*}
\end{thm}

\begin{prop} \label{t:approx}
Let $X$ be a K\"ahler manifold and let $\omega_0 \in H^0(X, \K_\infty^1)$ be the K\"ahler form, i.e., represented as a section by a choice of local K\"ahler potentials $\varphi_{0,m}$.  Then for any $\tau \in C^\infty(X)$ there exists $\sigma \in C^\omega(X)$ such that $\omega_0 + i \partial \dbar( \sigma - \tau)$ is still a positive form; equivalently, $\varphi_{0,m} + \sigma - \tau$ are strictly subpluriharmonic.
\end{prop}

\begin{proof}
Choose a $C^\omega$ embedding $X \hookrightarrow \Omega$ as above, take the open cover $\{ U_m \}$ as above and $\tilde{\eta} \in C^0(X)$ as in Lemma \ref{t:eta}.  Extend $\tilde{\eta}$ to a continuous positive function on $\Omega$ and extend $\tau$ to a smooth function $f \in C^\infty(\Omega)$.  Then by the Whitney approximation theorem, there exists $g \in C^\omega(\Omega)$ such that $h := f - g$ satisfies \eqref{C}.  Then setting $\sigma := g|_X \in C^\omega(X)$, we have for every $m \in \N$,
\begin{align*}
\varphi_{0,m} - h|_X = \varphi_{0,m} - (f|_X - g|_X) = \varphi_{0,m} - \tau + \sigma 
\end{align*}
is strictly plurisubharmonic.
\end{proof}

\begin{proof}[Proof of Proposition \ref{t:anrep}]
As before, we realise the K\"ahler form $\omega_0$ in terms of local K\"ahler potentials and hence think of it as a global section in $H^0(X, \K_\infty^1)$ with strictly plurisubharmonic representatives.

Taking the long exact sequences in cohomology associated to the diagram \eqref{E}, we get
\begin{align} \label{F}
\vcenter{
\xymatrix{
H^0(X, C_X^\omega) \ar[r] \ar[d] & H^0(X, \K_\omega^1) \ar[d]^s \ar[r]^-{\delta_\omega} & H^1(X, \O_X/i\R_X) \ar@{=}[d] \ar[r] & H^1(X, C_X^\omega) \ar[d] \\
H^0(X, C_X^\infty) \ar[r] & H^0(X, \K_\infty^1) \ar[r]_-{\delta_\infty} & H^1(X, \O_X/i\R_X) \ar[r] & H^1(X, C_X^\infty) .
} }
\end{align}
Now, $H^1(X, C_X^\infty) = 0$, since $C_X^\infty$ is a soft sheaf; it is also true, by a result of H.~Cartan, that $H^1(X, C_X^\omega) = 0$, see \cite[Th\'eor\`eme 3]{Cartan}. Hence, the maps $\delta_\infty$ and $\delta_\omega$ in \eqref{F} are both surjective.  Therefore, there exists $\omega^a \in H^0(X, \K_\omega^1)$ such that $\delta_\omega( \omega^a) = \delta_\infty( \omega_0)$, which in turn implies that there exists $\tau \in C^\infty(X)$ such that $\omega_0 - \omega^a = i \partial \dbar \tau$. Here, we are identifying $\omega^a$ with its image under $s$ in $H^0(X, \K_\infty^1)$. Proposition \ref{t:approx} then produces a $\sigma \in C^\omega(X)$ such that $\omega := \omega_0 + i \partial \dbar (\sigma - \tau) = \omega^a + i \partial \dbar \sigma$
is a positive, real analytic form. This concludes the proof.
\end{proof}

\vspace{0.4cm}
\begin{center}
 ----------------------------------------
\end{center}
\vspace{0.4cm}

\end{document}